\documentclass{nature}

\usepackage{epsfig,amssymb,amsmath,amsthm,amssymb,amsbsy,mathrsfs}
\usepackage{graphicx}
\usepackage[version=3]{mhchem} 
\usepackage{algorithm,algorithmic}
\usepackage{multirow}
\usepackage{xcolor}
\usepackage{lineno}
\usepackage{caption}
\usepackage{hyperref}

\allowdisplaybreaks[4]

\graphicspath{{fig/}}
\newenvironment{keyword}{Keywords:}{\par\medskip}
\newcommand*{\sep}{,}

\title{A fast low-rank inversion algorithm of dielectric matrix in GW approximation}
\author{Zhengbang Zhou\(^{1,}{}^{*}\),
Huanhuan Ma\(^2\),
Wentiao Wu\(^{3,4}\),
Weiguo Gao\(^{1,5,6}\),
Jinlong Yang\(^2\),
Meiyue Shao\(^{5,6,}{}^{*}\),
and Wei Hu\(^{2,3,4,}{}^{*}\)}

\newcommand{\code}[1]{\texttt{#1}}

\renewcommand{\vec}[1]{\mathbf{#1}}

\newcommand{\Ew}[1]{\mathrm{#1}}
\newcommand{\abs}[1]{\left\lvert#1\right\rvert}

\newcommand{\norm}[1]{\lVert#1\rVert}
\newcommand{\Fnorm}[1]{\lVert#1\rVert_{\mathrm{F}}}
\newcommand{\order}[1]{\mathcal{O}({#1})}
\newcommand{\occ}{\mathcal{V}}
\newcommand{\unocc}{\mathcal{C}}
\newcommand{\all}{\mathcal{N}}

\newcommand*{\md}{\mathop{}\mathopen{}\mathrm{d}}
\renewcommand{\Re}{\mathop{\mathrm{Re}}}

\newcommand*{\me}{\mathrm e}
\newcommand*{\mi}{\mathrm{i}}
\newcommand*{\R}{\mathbb{R}}
\newcommand*{\C}{\mathbb{C}}
\newcommand*{\GW}{\mathrm{GW}}

\newcommand*{\coc}{{C{_{\Ew{vc}}}\Omega^{-1}C_{\Ew{vc}}^{\mathsf{H}}}}
\DeclareMathOperator{\sgn}{sign}

\newcommand*{\trans}{^{\mathsf{T}}}
\newcommand*{\herm}{^{\mathsf{H}}}

\newcommand*{\conj}[1]{\overline{#1}}
\newcommand{\delete}[1]{}

\DeclareMathOperator{\diag}{diag}

\DeclareMathOperator{\sn}{sn}
\DeclareMathOperator{\cn}{cn}
\DeclareMathOperator{\dn}{dn}
\DeclareMathOperator{\tr}{tr}

\begin{document}

\maketitle

\begin{affiliations}
\item School of Mathematical Sciences, Fudan University, Shanghai 200433, China
\item Hefei National Research Center for Physical Sciences at the Microscale, University of Science and Technology of China, Hefei, Anhui 230026, China
\item School of Data Science, University of Science and Technology of China, Hefei, Anhui 230026, China
\item MOE Key Laboratory for Mathematical Foundations and Applications of Digital Technology, University of Science and Technology of China, Hefei, Anhui 230026, China
\item School of Data Science, Fudan University, Shanghai 200433, China
\item MOE Key Laboratory for Computational Physical Sciences, Fudan University, Shanghai 200433, China

\textit{E-mail: 
zbzhou21@m.fudan.edu.cn; myshao@fudan.edu.cn; whuustc@ustc.edu.cn}

\end{affiliations}

\begin{abstract}
The dielectric response function and its inverse are crucial physical quantities in materials science.
We propose an accurate and efficient strategy to invert the dielectric function matrix.
The GW approximation, a powerful approach to accurately describe many-body
excited states, is taken as an application to demonstrate accuracy and efficiency.
We incorporate the interpolative separable density fitting (ISDF) algorithm with Sherman--Morrison--Woodbury (SMW) formula to accelerate the inversion process by exploiting low-rank properties of dielectric function in plane-wave GW calculations.
Our ISDF--SMW strategy produces accurate quasiparticle energies with
\(O(N_{\Ew{r}}N_{\Ew{e}}^2)\) computational cost (\(N_{\Ew{e}}\) is the number
of electrons and \(N_{\Ew{r}}=100\)--\(1000N_{\Ew{e}}\) is the number of grid
points) with negligible small error of \(0.03\)~eV for both complex molecules and solids.
This new strategy for inverting the dielectric matrix
can be \(50\times\) faster than the current state-of-the-art implementation in
BerkeleyGW, resulting in two orders of magnitude speedup for total GW calculations.

\end{abstract}


\section{Introduction}
\label{sec:introduction}
The dielectric functions and polarizability functions are important physical
quantities in condensed matter physics and computational chemistry.
Experimentally, the dielectric function can be utilized to simulate optical
properties observed in experiments, such as absorption, reflection, and energy
loss~\cite{Liu17dielectric}.
On a theoretical level, many ground-state-based calculations concerning
excited-state necessitate the utilization of the dielectric
response function and polarizability function~\cite{Book:TDDFTSpringer,
SkoneSelfHybrid}.
The polarizability function is closely related to the dielectric function,
and its inversion through the Dyson equation and differential relation.
Computing the inverse of dielectric function matrix becomes one of the most
important steps in excited-state calculations.

There are two primary approaches to invert the dielectric matrix.
The conventional method involves explicitly constructing the irreducible
polarizability function, then the dielectric matrix, and finally inverting it.
The alternative approach involves obtaining the polarizability function
through the Sternheimer equation and subsequently deriving the dielectric
function~\cite{Giustino2010gw}.
Since the latter approach cannot explicitly obtain the dielectric function,
in most applications the first approach is preferred.
However, the dielectric function, \(\epsilon(\vec{r}, \vec{r}^{\prime})\),
becomes an \(N_{\Ew{r}}\times N_{\Ew{r}}\) matrix after discretization especially for plane waves.
Consequently, the conventional approach requires quartic time in terms of the
number of electrons \(N_{\Ew{e}}\) to construct the dielectric function
initially, followed by another cubic scaling time with a huge prefactor
(around \(10^9\)) to invert the dielectric function matrix.
This approach becomes very costly for large-scale calculations.
Since the inversion process is indispensable in most excited-state
calculations, the complexity of excited-state calculations becomes
prohibitively expensive for large systems.

\delete{Excited electrons play a vital role in determining crucial characteristics of
various materials, such as the efficiency of converting solar energy, protein
folding, photosynthesis reactions, and chemical 
reactions~\cite{He2019DFTforBattery, Makkar2021DFTnanomaterials,
Marzari2021ElecMaterial}.
Therefore, comprehending the properties of these electrons in their
excited-state such as the bandgap is essential for several advancing scientific fields.
However, accurate and efficient calculations to identify properties of electrons
comes with exponential time complexity with respect to the problem size.
The great time complexity renders the exact solution impractical, and some
approximations must be made.

First-principles computations use fundamental theories from quantum mechanics
to calculate wavefunctions and band properties.
Compared with other conventional approaches, first-principle calculations
introduce no correlations from experiments, 
and density-function theory (DFT) is one of the most practical
method~\cite{lin2019KSDFT}.
Nevertheless, as a ground-state theory,
DFT often seriously underestimate the bandgap of the
material~\cite{hybertsen1986electron}, which is crucial
in determining characteristics of various materials.
}

The GW method is an excited-state method rooted in many-body Green's function
theory.
It has become the gold standard for one-particle excitation energy
calculations due to its balance between efficiency and
accuracy~\cite{aryasetiawan1998gw, Hedin1965new, Onida2002Greens}.
Over the past decades, it has proven to be a powerful tool for describing
one-particle excitations in various systems, including
molecular~\cite{Kang2010Enhanced, ke2011all, liu2015numerical} and
solid-state materials~\cite{ma2007quasiparticle, spataru2004excitonic}.
Researchers have applied the GW method to investigate material science
hotspots such as low-dimensional nanostructures and
surfaces~\cite{Gaerosa2018GWapplication}.
Several methods and algorithms have been developed based on the GW method,
including G\(_{0}\)W\(_{0}\)~\cite{Hedin1965new},
GW\(\Gamma\)~\cite{Kuwahara2016GWGammaBSE},
scGW~\cite{schilfgaarde2006scGW}, among others.
G\(_0\)W\(_0\), as a variation of the GW method, simplifies the computation
by setting both the Green's function and the screened Coulomb interaction to
be one-shot, resulting in lower time complexity compared to other
alternatives.

Although G\(_0\)W\(_0\) calculation makes a concession to the accuracy,
it is still too expensive for many applications.
As an excited-state method, G\(_0\)W\(_0\) calculation also requires
inverting the dielectric function matrix.
This step introduces a quartic term and a cubic term with ultra high
factor in the complexity of the G\(_0\)W\(_0\) calculation.
BerkeleyGW, which was in the finalist of the ACM Gordon Bell Prize of
\(2020\)~\cite{ben2020GBBerkeleyGW} and is regarded the state-of-art
implementation of the GW calculation, also mentioned the difficulty introduced
by the huge prefactor.
For instance, it takes \(25.08\) seconds and \(162\)~MB RAM for the DFT
calculation with Quantum ESPRESSO~\cite{Giannozzi2009QE} to calculate all
bands of Si\(_{32}\) with a
single processor, while the static COHSEX G\(_0\)W\(_0\) calculation with
BerkeleyGW costs \(24454.35\) seconds and \(2519\)~MB RAM.

To reduce computational cost and memory usage, several low-scaling algorithms
have been proposed.
One such algorithm is the tensor hypercontraction (THC)
algorithm~\cite{lu2015compression},
known as the \emph{interpolative separable density fitting} (ISDF)
algorithm~\cite{hu2017interpolative, Lu2017CubicSA}.
The ISDF algorithm compresses orbital pair functions through a low-rank
approximation, significantly accelerating multi-center integral calculations
in quantum chemistry.
Moreover, the ISDF algorithm has been applied to develop an accurate and
efficient cubic scaling GW calculation strategy for both molecular and solid
systems, using atomic orbitals~\cite{gao2020accelerating} and a plane-wave
basis set~\cite{ma2021GWISDF}.
However, these implementations still require explicitly inverting the
dielectric function matrix, leaving the huge prefactor unresolved.
To address this issue, Liu et al.\ provide an approach to compute the inverse
using the Sherman--Morrison--Woodbury (SMW) formula~\cite{liu2015numerical}.
Although the prefactor when inverting the matrix is significantly reduced,
this process is only applicable to small systems that contain tens of atoms
because of its \(\order{N_{\Ew{e}}^6}\) time complexity and
\(\order{N_{\Ew{e}}^4}\) space complexity.

In this work, we present an accurate and effective cubic scaling algorithm by
combing the ISDF low-rank approximation and the SMW formula to invert the
dielectric function matrix \(\epsilon(\vec{r},\vec{r}^{\prime})\) with a small
prefactor and low memory usage.
We also apply this fast low-rank inversion algorithm to the G\(_0\)W\(_0\)
calculations as an application to excited-state calculations, illustrating the
accuracy and efficiency of our algorithm.
Table~\ref{table:notation} summarizes the notation of important variables and
operators in this work.


\section{Result and discussion}
\label{sec:numerical_experiments}
\subsection{Theoretical framework of fast inversion of dielectric matrix and
its application in GW calculation}

The dielectric function and polarizability functions are significant
physical quantity in computational chemistry and condensed matter physics,
both experimentally and theoretically.
The polarizability function is closely related to the inverse of
dielectric function.
Inverting the dielectric function matrix becomes one of the most important
steps in excited-state calculations.
The conventional process to construct the inverse of dielectric matrix is
first constructing dielectric matrix with quartic time complexity and cubic
space complexity, and then inverting it with cubic time complexity but a huge
prefactor.

We provide a cubic scaling method to invert the dielectric matrix accurately
and efficiently with quadratic space complexity and a modest prefactor.
We use the \emph{ISDF algorithm} to generate polarizability operator and
dielectric matrix with low-rank properties, so that the \emph{SMW formula} can
be adopted to invert the matrix.
To eventually achieve cubic scaling, we apply \emph{Cauchy integral} to the
coupling part of the dielectric matrix.
A detailed description of our fast low-rank inversion algorithm can be found
in Section~\ref{sec:algorithm}.

\begin{figure}[!tb]
\centering
\includegraphics[width=0.90\textwidth]{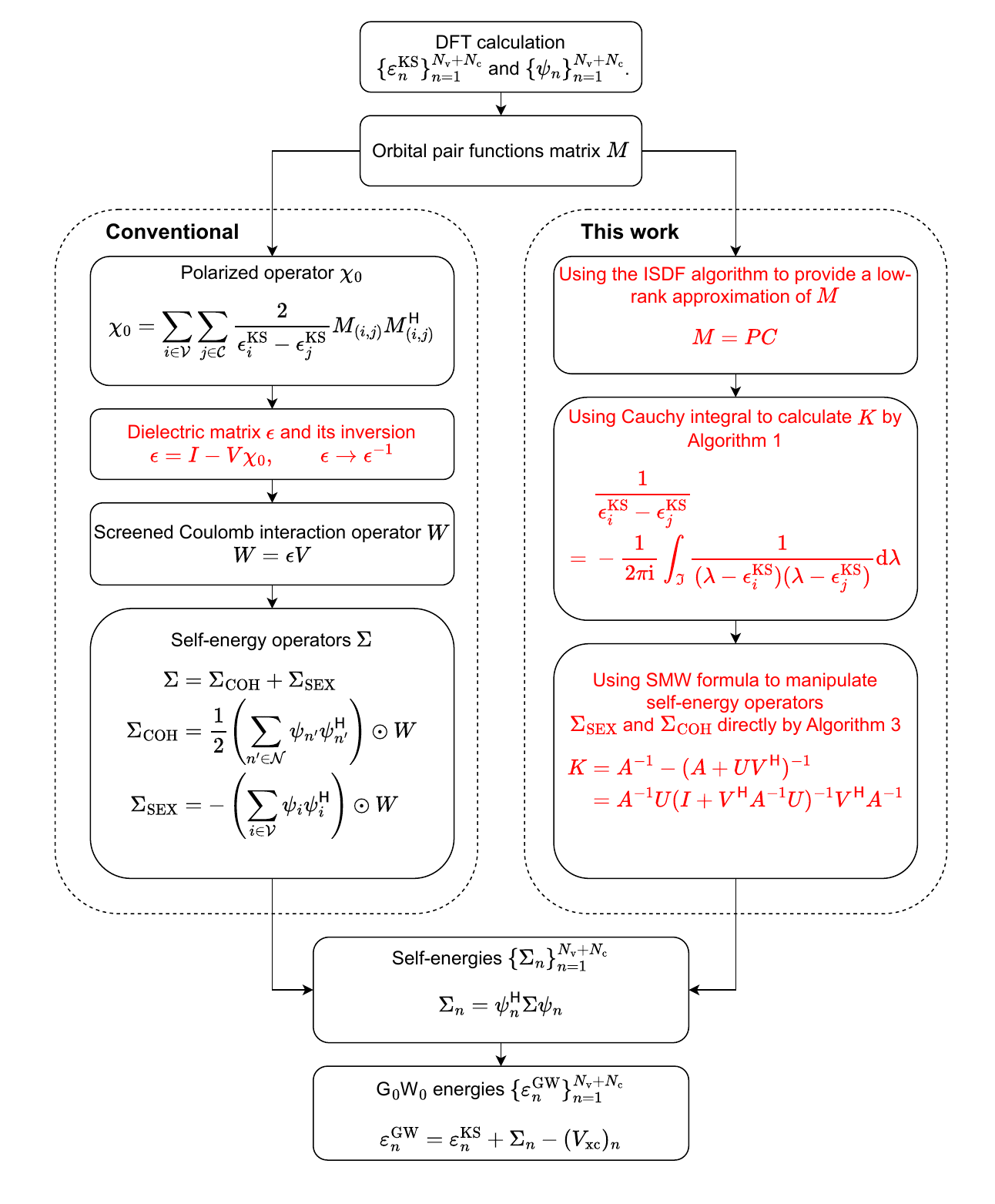}
\caption{Flowchart difference between conventional and low-rank G\(_0\)W\(_0\)
calculations.
The inversion step is highlighted.}\label{fig:GW_flowchart}
\end{figure}

The GW method is a well-developed method to calculate the excited-state
energies.
The framework of GW method introduces self-energies to go beyond the
mean-field approximation in DFT calculation and approximate the self-energies
from a set of Hedin's equations, which is used to express screened Coulomb
effects.
Inverting the dielectric matrix, however, is also one of the most
time-consuming and memory-intensive step in the GW method.
The conventional approach involves at least quartic time complexity and cubic
space complexity.
We apply our fast inversion algorithm to the GW method, reducing one-order of
both the time complexity and the space complexity compared to the conventional
approach.
An overview of our improved GW strategy and the comparison to the conventional
one is shown in Fig.~\ref{fig:GW_flowchart}.
Detailed descriptions and theoretical derivation on the stability of our
improved GW calculation strategy can be found in Section~\ref{sec:algorithm}.

In the numerical examples, we use both molecular systems (C\(_{60}\)) and
solid systems (bulk silicon and SrTiO\(_3\)) under the periodic boundary
condition to illustrate the accuracy, efficiency, and generally applicability
of our strategy.
All reported results except examples arguing complexity issues are performed
on a single node with \(16\) computational cores of a \(800\)~MHz processor
with \(4\)~TB memory using MATLAB, while the examples showing complexity are
limited on a single core.
All experiments are based on the solution provided by
\code{KSSOLV}~\cite{Jiao2022KSSOLV2.0} for the
corresponding KS--DFT problem.
To highlight the scaling and minimize the impact of less important factors, in
the experiments we set the number of the unoccupied orbitals \(N_{\Ew{c}}\) to
be equal to the number of the occupied orbitals \(N_{\Ew{v}}\).
Parameters of our test systems are given in Table~1 of
Supplementary information.

\subsection{Errors of self-energies}
\label{subsec:SelfEnergyError}
Our improved method produces approximation errors of self-energies from both
the ISDF algorithm and numerical quadrature of Cauchy integral, and are
bounded with
\begin{equation}
\begin{aligned}
\mathcal{O}\bigl(\delta M + \me^{-N_{\lambda}/C}\bigr).
\label{eq:rel_err_result}
\end{aligned}
\end{equation}
Here, \(\delta M\) refers to the errors introducing by the ISDF algorithm to
the orbital pair functions matrix~\(M\), \(N_{\lambda}\) is the number of
quadrature nodes in Cauchy integral, and \(C\) is the condition number related
to the bandgap of the system.
More precise expressions of these terms in the error bound will be given
in~\eqref{eq:rel_err_Cauchy} and~\eqref{eq:rel_err_ISDF}, respectively,
and the error analysis can be found in Section~\ref{subsec:errorAnalysis}.
We highlight two key facts here.
First, an accurate strategy in the ISDF algorithm is necessary because the
error bound in equation~\eqref{eq:rel_err_result} requires a good
approximation of \(M\).
Second, the bound given by~\eqref{eq:rel_err_result} is pessimistic.
Numerical experiments show that a much accurate result can be expected for
real systems.

We first investigate the effect of errors in the ISDF algorithm.
Our improved strategy applies the ISDF algorithm to the orbital pair functions
matrix \(M\) of different sets.
Different combinations of \(M\) are used in different steps in our improved
strategy.
The product of occupied and unoccupied states is used to construct operators
in the fast-inversion strategy, and other combinations are used to construct
four-center two-electron integrals and manipulate self-energies in our
improved GW strategy.
We need to treat them separately since they have different effects on the
errors of self-energies.
To illustrate the effect of errors in the ISDF algorithm, we vary the ISDF
coefficients in different combinations of \(M\).
We first change the ISDF coefficient with respect to occupied and unoccupied
states, denoted by \(k_{\Ew{vc}}\), and fix the other ISDF coefficients,
denoted by \(k_{\Ew{vn}}\) and \(k_{\Ew{nn}}\) to~\(8.0\).
Then we fix the ISDF coefficient of the first one to \(8.0\) and vary the
other ISDF coefficients.

\begin{figure}[!tb]
\centering
\includegraphics[width=0.97\textwidth]{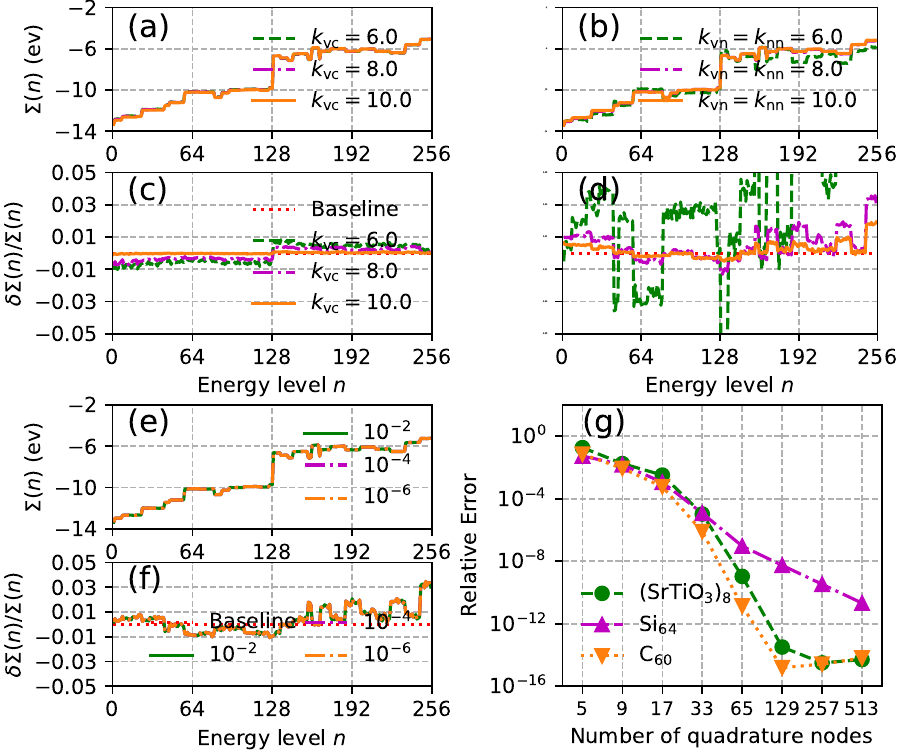}
\caption{
Errors of self-energies for Si\(_{64}\) in improved GW calculation strategy.
(a)~and~(c): The ISDF coefficient \(k_{\Ew{vn}}\) and \(k_{\Ew{nn}}\)
are fixed to \(8.0\), and \(k_{\Ew{vc}}\) varies from \(6.0\) to \(10.0\).
(b)~and~(d): The ISDF coefficient \(k_{\Ew{vc}} = 8.0\) is fixed,
and \(k_{\Ew{vn}}=k_{\Ew{nn}}\) varies from \(6.0\) to \(10.0\).
The green, purple, and orange lines in~(d) represent
\(k_{\Ew{vn}} = k_{\Ew{nn}} = 6.0\), \(k_{\Ew{vn}} = k_{\Ew{nn}} = 8.0\),
and \(k_{\Ew{vn}} = k_{\Ew{nn}} = 10.0\), respectively.
(a)~to~(d) do not involve Cauchy integral. 
(e)~and~(f): The ISDF coefficients
\(k_{\Ew{vc}} = k_{\Ew{vn}} = k_{\Ew{nn}} = 8.0\) are
fixed, and the upper bound of the relative error of Cauchy integral result
in Algorithm~\ref{algo:COmegaC} varies from \(10^{-2}\) to \(10^{-6}\).
(g): Errors of Cauchy integral in terms of the number of quadrature nodes.
}\label{fig:ISDF_k_vc}
\end{figure}

The results, as shown in Figs.~\ref{fig:ISDF_k_vc}~(a) to (d), reveal the
following two facts.
First, increasing any ISDF coefficients can reduce the relative errors on
self-energies as expected.
Second, comparing Figs.~\ref{fig:ISDF_k_vc}~(c) and~(d) when \(k = 6.0\),
we notice that the taking \(k_{\Ew{vn}} = k_{\Ew{nn}} = 6.0\) produces much
larger relative error compared to taking \(k_{\Ew{vc}} = 6.0\).
Hence, self-energies are much more sensitive to \(k_{\Ew{vn}}\) and
\(k_{\Ew{nn}}\) than \(k_{\Ew{vc}}\).
The reason is that \(k_{\Ew{vc}}\) affects the operators' calculation, which
has only indirect impact to self-energies, whereas \(k_{\Ew{vn}}\) and
\(k_{\Ew{nn}}\) affect the self-energies calculation straightway.
It is worth mentioning that both time complexity and space complexity
of our method is quadratic in terms of the ISDF coefficients (see
Table~\ref{table:TimeSpacecomplexity}).
Balancing is needed between accuracy and efficiency.
The experiment above shows that it is possible to choose a smaller value of
\(k_{\Ew{vc}}\) and larger values of \(k_{\Ew{vn}}\) and \(k_{\Ew{nn}}\) to
accelerate the calculation while ensuring the relative accuracy of the
self-energies.


Then we show the relationship between the quadrature error in Cauchy integral
and the error of self-energies.
Algorithm~\ref{algo:COmegaC} explains how to adopt Cauchy integral in our
strategy with adaptive refinement.
The algorithm provides a result based on a given threshold that also bounds
the quadrature error.
We still use Si\(_{64}\) as a case study and set different upper bounds
on the relative error in Algorithm~\ref{algo:COmegaC} to demonstrate the
relationship between the errors.
We present our computational results in Figs.~\ref{fig:ISDF_k_vc}~(e)~to~(g).
Figs.~\ref{fig:ISDF_k_vc}~(e) and~(f) show that with quadrature error
ranging from \(10^{-2}\) to \(10^{-6}\), the corresponding error curves for
the self-energies are closely intertwined.
This observation suggests that the error of Cauchy integral has marginal
impact on the accuracy of the self-energies.
This allows us to choose a loose error estimate for
Algorithm~\ref{algo:COmegaC} without impacting the accuracy of self-energies.

%

Our quadrature rule is an adaptive one that automatically refines the
partition.
As the quadrature error decreases geometrically with respect to the number of
quadrature nodes, the cost of numerical quadrature is moderately small.
Numerical tests in Fig.~\ref{fig:ISDF_k_vc}~(g) agree with our
analysis.
By quadaturing with \(513\) nodes,
the quadrature error already drops to a
satisfactory level, indicating that it is easy to achieve a negligible
relative error of Cauchy integral for both molecular and solid systems.
Hence it is not really important to balance between efficiency and accuracy
in calculating Cauchy integral.

\subsection{Low-rank properties for functions and operators}
\label{subsec:lowRank}
The ISDF algorithm and Sherman--Morrison--Woodbury formula assume that
both the orbital pair functions matrix \(M\) and the polarizability matrix
\(\chi\) have certain low-rank properties.
This assumption is usually plausible in practice.
We take Si\(_{64}\) and (SrTiO\(_3\))\(_8\) as examples to demonstrate the
decay of singular values for both \(M\) and~\(\chi\).

\begin{figure}[!tb]
\centering
\includegraphics[width=0.99\textwidth]{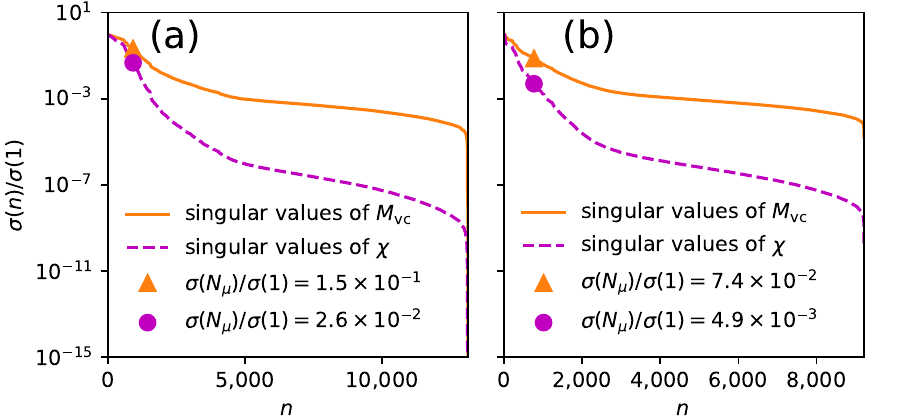}
\caption{Singular values of conventional \(M_{\Ew{vc}}\) and \(\chi\) in
(a) Si\(_{64}\) and (b) {(SrTiO\(_3\))}\(_8\).
Here, \(k_{\Ew{vc}}=8.0\), \(N_{\mu}=k_{\Ew{vc}}\sqrt{N_{\Ew{v}}N_{\Ew{c}}}\),
\(N_{\Ew{v}}\) and \(N_{\Ew{c}}\) are given in
Table~1 of Supplementary information.}\label{fig:ISDF_chi_svd}
\end{figure}

Numerical experiments in Fig.~\ref{fig:ISDF_chi_svd} show that when
\(k_{\Ew{vc}} = 8.0\), the ISDF algorithm
introduces a relative error of \(\Theta(10^{-1})\) to \(M_{\Ew{vc}}\) in both
systems and a relative error of \(\Theta(10^{-2})\) to \(\chi\) in Si\(_{64}\).
The decay pattern of the singular values is system-dependent, as illustrated
in Fig.~\ref{fig:ISDF_chi_svd}.
The singular values of both \(\chi\) and \(M_{\Ew{vc}}\) in Si\(_{64}\) decay
rapidly but smoothly, while \(\chi\) in (SrTiO\(_3\))\(_8\) has a much faster
decay in its leading singular values.
This rapid-decay phenomenon is strongly associated with the pattern of orbital
pair functions matrix \(M_{\Ew{vc}}\), as the leading singular values of the
orbital pair functions matrix \(M_{\Ew{vc}}\) of (SrTiO\(_3\))\(_8\) decrease
much faster than that of Si\(_{64}\).

We remark that though errors introduced by the ISDF algorithm to \(M\) and
\(\chi\) can be pretty large at the first glance, results in
Fig.~\ref{fig:ISDF_k_vc} suggest that the self-energies can still be computed
with satisfactory accuracy.

\subsection{Computational complexity}
\label{subsec:Complexity}
Silicon bulk systems are known for their symmetry and highly extendible
properties.
In this study, we demonstrate the effectiveness of our method using Si\(_8\),
Si\(_{16}\), Si\(_{32}\), and Si\(_{64}\), in terms of time complexity.
We compare our strategy with the conventional G\(_0\)W\(_0\) calculation
strategy, the strategy by Ma et al.~\cite{ma2021GWISDF}, all implemented in
MATLAB.
The test result of the Fortran-based package BerkeleyGW, which uses the
conventional G\(_0\)W\(_0\) strategy, is also listed for reference.
In this test, we only use a single processor.

\begin{figure}[!tb]
\centering
\includegraphics[width=0.99\textwidth]{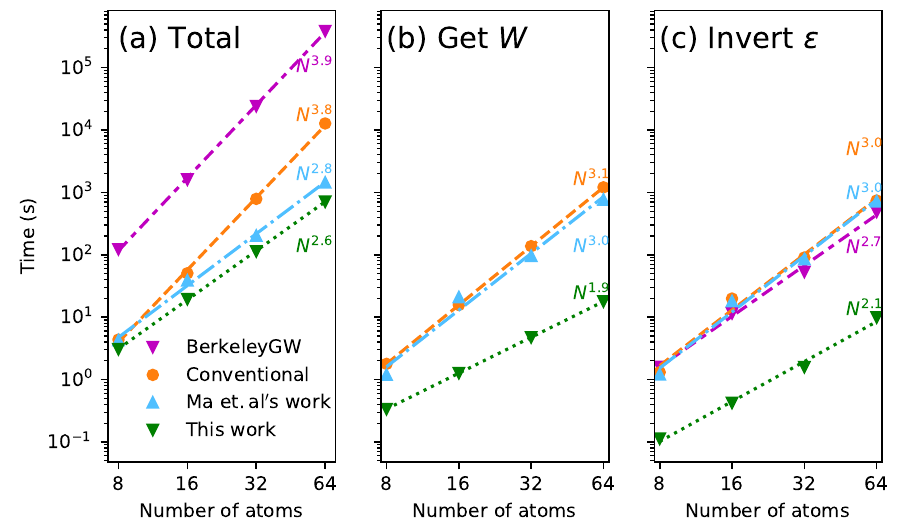}
\caption{
Total time for the conventional strategy, strategy from~\cite{ma2021GWISDF}
and our improved strategy on Si with different bulk sizes.
The ISDF coefficients are set to be \(8.0\), and the threshold for Cauchy
integral is set to be \(10^{-7}\).
(a) Execution time from calculating all operators to calculating self-energies
in G\(_0\)W\(_0\) calculation.
The time for ISDF algorithm is included.
(b) Execution time until obtaining the screened Coulomb interaction matrix $W$,
excluding the time for ISDF.
BerkeleyGW is not listed here as its implementation is a bit different so that
it is not easy to make a fair comparison.
(c) Total execution time for inverting \(\epsilon\).
}\label{fig:ISDF_Si_bulk}
\end{figure}

Fig.~\ref{fig:ISDF_Si_bulk} shows that our strategy significantly reduces
the time required for G\(_{0}\)W\(_{0}\) calculation compared to other
strategies.
Our main improvement in the inversion of the dielectric matrix leads to about
\(10\times\) speedup compared to other alternatives already for these small
systems.
The benefit of our strategy becomes more and more significant as the problem
size grows, because our strategy also has a lower asymptotic scaling than
other strategies.

We list the complexity of each algorithm in
Table~\ref{table:TimeSpacecomplexity}
(BerkeleyGW is as same as the conventional strategy).
Let us take the strategy from~\cite{ma2021GWISDF} as an example.
The time complexity of this strategy is
\(\order{N_{\Ew{r}}^3 + k_{\Ew{vc}}^2N_{\Ew{e}}^4}\), where
\(N_{\Ew{r}} \sim 1000 N_{\Ew{e}}\).
The cubic term dominates for small systems with tens of electrons
because \(N_{\Ew{r}}\) is very large compared to \(N_{\Ew{e}}\).
The corresponding curve in Fig.~\ref{fig:ISDF_Si_bulk}~(a) stays low, and is
expected to grow rapidly for larger systems.

We also notice that BerkeleyGW, a Fortran-based high-performance computing
package, performs much more slowly than our MATLAB-based version.
The main reason is that several computationally intensive steps in BerkeleyGW
do not utilize level-3 BLAS, resulting in a significant loss of computational
efficiency.
As a well-developed package, BerkeleyGW supports multiple \(k\)-points
calculation, and manages to use a unified approach to handle the
\(\Gamma\)-point (also known as the single \(k\)-point) and multiple
\(k\)-points.
Such a unified approach is based on several do loops, making the most
computational intensity part of the self-energies calculation elementwise.
However, in the MATLAB version we do not support \(k\)-points yet.
This allows us to easily perform the computation using well-developed level-3
BLAS.

\section{Methods}
\label{sec:algorithm}
\subsection{GW approximation}
GW approximation expresses one-particle excitation energies with
quasiparticle equation as
\begin{equation}
\begin{aligned}
\Bigl(-\frac12\nabla^2 + V_{\Ew{ion}}(\vec{r}) + V_{\Ew{H}}(\vec{r})\Bigr)
\psi_i^{\Ew{QP}}(\vec{r})
+ \int_{\R^3} \md\vec{r}^{\prime}
\Sigma\bigl(\vec{r}, \vec{r}^{\prime}; \varepsilon_i^{\Ew{GW}}\bigr)
\psi_{i}^{\Ew{QP}}(\vec{r}^{\prime})
= \varepsilon_{i}^{\Ew{GW}} \psi_{i}^{\Ew{QP}}(\vec{r}),
\label{eq:quasiparticle}
\end{aligned}
\end{equation}
where \(\Sigma\) is the self-energy operator in frequency
space~\cite{Hedin1965new,Hedin1970effects}, \(\psi^{\Ew{QP}}_i\)'s are the
quasiparticle wavefunctions,
and \(\varepsilon_i^{\Ew{GW}}\)'s refer to quasiparticle energies.
We assume spin deformation and zero-temperature in this work so that the spin
variables are ignored.
Also, we only consider the case of \(\Gamma\)-point, and hence indices of
\(k\)-points are excluded.

The self-energy operator can be expressed by the Green's functions~\(G\)
and the screened Coulomb potential~\(W\) as
\[\begin{aligned}
\Sigma(\vec{r}, \vec{r}^{\prime}; \omega)
={}& \frac{\mi}{2\pi} \int_{\R} \mathrm{e}^{-\mi\omega^{\prime}\eta}
G(\vec{r}, \vec{r}^{\prime};\omega-\omega^{\prime})
W(\vec{r},\vec{r}^{\prime}; \omega^{\prime}) \md \omega^{\prime}.
\end{aligned}\]
Here, \(\eta\) is a small frequency shift to guarantee convergence.
We use Lehmann representation of the one-particle Green's functions~\(G\)
with a set of quasiparticle wavefunctions
\begin{equation}\begin{aligned}
G(\vec{r}, \vec{r}^{\prime}; \omega)
= \sum_{n\in\all}
\frac{\psi_i^{\Ew{QP}}(\vec{r}) \conj{\psi_i^{\Ew{QP}}(\vec{r}^{\prime})}}
{\omega - \varepsilon_i^{\Ew{QP}} + \mi\sgn(\varepsilon_n - \mu)\eta},
\end{aligned}\end{equation}
and use a set of Hedin's equations to generate the screened Coulomb potential
operator~\(W\) as
\begin{equation}
\label{eq:Hedin_fourier}
\begin{split}
W(\vec{r}, \vec{r}^{\prime}; \omega)
={}& \int_{\R^3}\md \vec{r}^{\prime\prime}
\epsilon^{-1}(\vec{r}, \vec{r}^{\prime\prime} ; \omega)
V(\vec{r}^{\prime\prime}, \vec{r}^{\prime}) , \\
\epsilon(\vec{r}, \vec{r}^{\prime}; \omega)
={}& \delta(\vec{r} - \vec{r}^{\prime})
- \int_{\R^3}\md\vec{r}^{\prime\prime} V(\vec{r}, \vec{r}^{\prime\prime})
\chi_0(\vec{r}^{\prime\prime},\vec{r}^{\prime}; \omega) , \\
\chi_0(\vec{r}, \vec{r}^{\prime}; \omega)
={}& \sum_{i\in\occ} \sum_{j\in\unocc}
\psi^{\Ew{QP}}_i(\vec{r}) \conj{\psi^{\Ew{QP}}_j(\vec{r})}
\conj{\psi^{\Ew{QP}}_i(\vec{r}^{\prime})} \psi^{\Ew{QP}}_j(\vec{r}^{\prime})
\cdot \frac{1}{\omega + \varepsilon_i - \varepsilon_j + \mi\eta} \\
&- \sum_{i\in\occ}\sum_{j\in\unocc}
\conj{\psi^{\Ew{QP}}_i(\vec{r})} {\psi^{\Ew{QP}}_j(\vec{r})}
{\psi^{\Ew{QP}}_i(\vec{r}^{\prime})} \conj{\psi^{\Ew{QP}}_j(\vec{r}^{\prime})}
\cdot \frac{1}{\omega - \varepsilon_i + \varepsilon_j - \mi\eta}.
\end{split}
\end{equation}
The operators on the left-hand side are, respectively, the screened Coulomb
interaction operator, the dielectric operator, and the polarizability operator,
and \(V(\vec{r}, \vec{r}^{\prime})\) refers to the Coulomb interaction.

The simplest GW calculation is the one-shot GW, also known as G\(_0\)W\(_0\).
The standard procedure of G\(_0\)W\(_0\) is first getting the ground-states
information from DFT, and then replacing quasiparticle energies and
wavefunctions in equations above with the ground-states ones.
G\(_0\)W\(_0\) can be implemented by several commonly used approximations,
including the static COulomb Hole plus Screened EXchange (COHSEX)
approximation~\cite{Hedin1965new, Hedin1970effects},
the full-frequency
approximation~\cite{Godby1988GWfullfreqcd, Rojas1995GWfullfreqra},
and the generalized plasmon--pole approximation
(GPP)~\cite{hybertsen1986electron}.
To show the feasibility of our low-rank strategy of inverting the dielectric
matrix \(\epsilon\) in the framework of GW calculation, we use the static
COHSEX approximation.
In the static COHSEX approximation, we set the frequency to be \(0\) in the
screened Coulomb interaction operator \(W\) and let the Green's function \(G\)
to be a non-interacting one.

The self-energy operator \(\Sigma\) under the static COHSEX approximation can
be separated into
\[
\Sigma_{\Ew{COHSEX}}(\vec{r}, \vec{r}^{\prime})
=\Sigma_{\Ew{COH}}(\vec{r}, \vec{r}^{\prime})
+\Sigma_{\Ew{SEX}}(\vec{r}, \vec{r}^{\prime}),
\]
where the screened exchange part (SEX) and Coulomb hole part (COH),
respectively, are expressed as
\begin{equation}
\label{eq:Sigma_G0W0}
\begin{aligned}
\Sigma_{\Ew{COH}}(\vec{r}, \vec{r}^{\prime})
&= \frac{1}{2}\delta(\vec{r}-\vec{r}^{\prime})
\bigl(W(\vec{r}, \vec{r}^{\prime};\omega=0)
- V(\vec{r},\vec{r}^{\prime})\bigr),\\
\Sigma_{\Ew{SEX}}(\vec{r}, \vec{r}^{\prime})
&= -\sum_{i \in \occ} \psi_{i}(\vec{r}) \conj{\psi_{i}(\vec{r}^{\prime})}
W(\vec{r}, \vec{r}^{\prime};\omega=0).
\end{aligned}
\end{equation}
{The bare Coulomb part can be seperated from the screened exchange part of the
self-energy operator \(\Sigma_{\Ew{SEX}} = \Sigma_{\Ew{SEX\_X}} +
\Sigma_{\Ew{X}}\), and the corresponding self-energies of the band level
\(n\) for different parts are
\begin{equation}
\label{eq:AllPartSelfEnergies}
\begin{aligned}
\Sigma_{\Ew{SEX\_X}}(n)
={}& -\int_{\R^3\times\R^3} \md\vec{r}\md\vec{r}^{\prime}
\Bigl(\sum_{i\in{\occ}} \psi_{i}(\vec{r})\conj{\psi_{i}(\vec{r}^{\prime})}
\Bigr)
\psi_n(\vec{r}^{\prime})\conj{\psi_n(\vec{r})}
\bigl(W(\vec{r}, \vec{r}^{\prime}) - V(\vec{r}, \vec{r}^{\prime})\bigr),\\
\Sigma_{\Ew{X}}(n)
={}& -\int_{\R^3\times\R^3} \md\vec{r}\md\vec{r}^{\prime}
\Bigl(\sum_{i\in{\occ}} \psi_{i}(\vec{r})\conj{\psi_{i}(\vec{r}^{\prime})}
\Bigr)
\psi_n(\vec{r}^{\prime})\conj{\psi_n(\vec{r})}
V(\vec{r}, \vec{r}^{\prime}),\\
\Sigma_{\Ew{COH}}(n)
={}& \frac12\int_{\R^3\times\R^3} \md\vec{r}\md\vec{r}^{\prime}
\Bigl(\sum_{i\in{\all}} \psi_{i}(\vec{r})\conj{\psi_{i}(\vec{r}^{\prime})}
\Bigr)
\psi_n(\vec{r}^{\prime})\conj{\psi_n(\vec{r})}
\bigl(W(\vec{r}, \vec{r}^{\prime}) - V(\vec{r}, \vec{r}^{\prime})\bigr).\\
\end{aligned}
\end{equation}
}

For simplicity, we assume time-reversal symmetry in the system.
For the calculation with only \(\Gamma\)-point, time-reversal
symmetry becomes
\begin{equation}
\label{eq:timeReversalSymmetric}
\conj{\psi_n(\vec{r})} = \psi_n(\vec{r}),\quad \forall n,
\end{equation}
which means the wavefunctions are real in real space.
However, we continue to treat the wavefunctions \(\psi(\vec{r})\) as a complex
value to align with the general G\(_0\)W\(_0\) framework.
\delete{
Under the time-reversal symmetry assumption~\eqref{eq:timeReversalSymmetric},
we obtain the matrix forms of all operators in~\eqref{eq:Hedin_G0W0} as
\begin{equation}
\label{eq:matrices}
\chi_0 = 2{M_{\Ew{vc}}} \Omega^{-1} ({M_{\Ew{vc}}})\herm,\qquad
\epsilon = I - V\chi_0,\qquad
W = \epsilon^{-1}V,\\
\end{equation}
and components of the self-energies in~\eqref{eq:AllPartSelfEnergies} become
\[\begin{aligned}
\Sigma_{\Ew{SEX\_X}}(n)
&= -\sum_{i\in\occ}
\bigl(\psi_n \odot \conj{\psi_i}\bigr)\herm
(W - V) \bigl(\psi_n \odot \conj{\psi_i}\bigr),\\
\Sigma_{\Ew{X}}(n)
&= -\sum_{i\in\occ}
\bigl(\psi_n \odot \conj{\psi_i}\bigr)\herm V
\bigl(\psi_n \odot \conj{\psi_i}\bigr),\\
\Sigma_{\Ew{COH}}(n)
&= \frac{1}{2}\sum_{n^{\prime}\in\all}
\bigl(\psi_n \odot \conj{\psi_{n^{\prime}}}\bigr)\herm (W-V)
\bigl(\psi_n \odot \conj{\psi_{n^{\prime}}}\bigr).
\end{aligned}\]
In~\eqref{eq:matrices}, \(\Omega\) refers to a diagonal matrix with elements
\begin{equation}\begin{aligned}
\label{eq:Omega}
\Omega_{(i,j; k,l)} = \delta_{i,k}\delta_{j,l}
(\varepsilon_i^{\Ew{KS}} - \varepsilon_j^{\Ew{KS}}),
\qquad i,k\in\occ,\quad j,l\in\unocc,
\end{aligned}\end{equation}
and the matrix \(M_{\Ew{vc}}\) refers to orbital pair functions of all
occupied and unoccupied states, i.e.,
\[
\bigl(M_{\Ew{vc}}\bigr)_{(n; i,j)}
= \conj{\psi_i(\vec{r}_n)} {\psi_j(\vec{r}_n)},
\qquad i\in\occ,\quad j\in\unocc.
\]
Here, the pair \((i,j)\) together forms the column index of the matrix.
}

We compute the integral in real space using numerical quadrature based
on a set of uniformly distributed nodes in real space.
Numerical quadrature introduces summation over real-space indices and some
coefficients.
For the sake of simplicity in following paragraphs, we drop these coefficients.
We replace the integral over real space by numerical quadrature.
Applying numerical quadrature in real space to the operators
in~\eqref{eq:Hedin_fourier} under COHSEX approximation and the time-reversal
symmetry assumption~\eqref{eq:timeReversalSymmetric} leads to
\begin{equation}
\label{eq:Hedin_G0W0}
\begin{aligned}
\chi_0(\vec{r}, \vec{r}^{\prime})
={}& 2\sum_{i \in \occ} \sum_{j \in \unocc}
\psi_i(\vec{r}) \conj{\psi_i(\vec{r}^{\prime})}
\psi_j(\vec{r}^{\prime})\conj{\psi_j(\vec{r})}
\cdot\frac{1}{\varepsilon_i^{\Ew{KS}} - \varepsilon_j^{\Ew{KS}}},\\
\epsilon(\vec{r}, \vec{r}^{\prime})
={}& \delta(\vec{r}-\vec{r}^{\prime})
- \sum_{\vec{r}^{\prime\prime}} V(\vec{r}, \vec{r}^{\prime\prime})
\chi_0(\vec{r}^{\prime\prime}, \vec{r}^{\prime}),\\
W(\vec{r}, \vec{r}^{\prime})
={}& \sum_{\vec{r}^{\prime\prime}}
\epsilon^{-1}(\vec{r}, \vec{r}^{\prime\prime})
V(\vec{r}^{\prime\prime}, \vec{r}^{\prime}),
\end{aligned}
\end{equation}
and the components of the self-energies in~\eqref{eq:AllPartSelfEnergies} become
\begin{equation}\begin{aligned}
\Sigma_{\Ew{SEX\_X}} (n)
={}& -\sum_{i\in\occ} \sum_{\vec{r},\vec{r}^{\prime}}
\sum_{\vec{r}^{\prime\prime}}
\psi_i(\vec{r}) \conj{\psi_n(\vec{r})}
\conj{\psi_i(\vec{r}^{\prime})}\psi_n(\vec{r}^{\prime})
\bigl(\epsilon^{-1}(\vec{r}, \vec{r}^{\prime\prime})
- \delta(\vec{r}, \vec{r}^{\prime\prime}) \bigr)
V(\vec{r}^{\prime\prime}, \vec{r}^{\prime}),\\
\Sigma_{\Ew{X}} (n)
={}& -\sum_{i\in\occ} \sum_{\vec{r},\vec{r}^{\prime}}
\psi_i(\vec{r}) \conj{\psi_n(\vec{r})}
\conj{\psi_i(\vec{r}^{\prime})}\psi_n(\vec{r}^{\prime})
V(\vec{r}^{\prime}, \vec{r}^{\prime}),\\
\Sigma_{\Ew{COH}}(n)
={}&\frac{1}{2} \sum_{n'\in{\all}} \sum_{\vec{r}, \vec{r}^{\prime}}
\sum_{\vec{r}^{\prime\prime}}
\psi_{n'}(\vec{r})\conj{\psi_n(\vec{r})}
\conj{\psi_{n'}(\vec{r}^{\prime})}\psi_n(\vec{r}^{\prime})
\bigl(\epsilon^{-1}(\vec{r}, \vec{r}^{\prime\prime} )
- \delta(\vec{r}, \vec{r}^{\prime\prime})\bigr)
V(\vec{r}^{\prime\prime}, \vec{r}^{\prime}),
\label{eq:selfenergySum}
\end{aligned}\end{equation}
where \(\vec{r}\), \(\vec{r}^{\prime}\), and \(\vec{r}^{\prime\prime}\) are
taken to represent points on the set of real-space grids.

We prefer the matrix forms of all operators, where the low-rank properties
can be more intuitively observed.
We can obtain the matrix forms of operators and self-energies by simply
replacing summation in~\eqref{eq:Hedin_G0W0} and~\eqref{eq:selfenergySum} as
\begin{equation}\begin{aligned}
\label{eq:matrices}
\chi_0 = 2{M_{\Ew{vc}}} \Omega^{-1} ({M_{\Ew{vc}}})\herm,\qquad
\epsilon = I - V\chi_0,\qquad
W = \epsilon^{-1}V,\\
\end{aligned}\end{equation}
\begin{equation*}\begin{aligned}
\Sigma_{\Ew{SEX\_X}}(n)
&= -\sum_{i\in\occ}
\bigl(\psi_n \odot \conj{\psi_i}\bigr)\herm
\bigl((\epsilon^{-1} - I) V\bigr)
\bigl(\psi_n \odot \conj{\psi_i}\bigr),\\
\Sigma_{\Ew{X}}(n)
&= -\sum_{i\in\occ}
\bigl(\psi_n \odot \conj{\psi_i}\bigr)\herm V
\bigl(\psi_n \odot \conj{\psi_i}\bigr),\\
\Sigma_{\Ew{COH}}(n)
&= \frac{1}{2}\sum_{n^{\prime}\in\all}
\bigl(\psi_n \odot \conj{\psi_{n^{\prime}}}\bigr)\herm 
\bigl((\epsilon^{-1} - I) V\bigr)
\bigl(\psi_n \odot \conj{\psi_{n^{\prime}}}\bigr).
\end{aligned}\end{equation*}
In~\eqref{eq:matrices}, \(\Omega\) refers to a diagonal matrix with elements
\begin{equation}\begin{aligned}
\label{eq:Omega}
\Omega_{(i,j; k,l)} = \delta_{i,k}\delta_{j,l}
(\varepsilon_i^{\Ew{KS}} - \varepsilon_j^{\Ew{KS}}),
\qquad i,k\in\occ,\quad j,l\in\unocc,
\end{aligned}\end{equation}
and the matrix \(M_{\Ew{vc}}\) refers to orbital pair functions of all
occupied and unoccupied states, i.e.,
\[
\bigl(M_{\Ew{vc}}\bigr)_{(n; i,j)}
= \conj{\psi_i(\vec{r}_n)} {\psi_j(\vec{r}_n)},
\qquad i\in\occ,\quad j\in\unocc.
\]
Here, the pair \((i,j)\) together forms the column index of the matrix.

As we have mentioned before, both time complexity and space complexity of the
G\(_0\)W\(_0\) calculation are pretty high.
Evaluating \(\chi_0\) as a matrix costs \(\order{N_{\Ew{e}}^4}\);
inverting \(\epsilon\) costs \(\order{N_{\Ew{e}}^3}\) with a huge prefactor;
and calculating \(\Sigma\) conventionally costs \(\order{N_{\Ew{e}}^5}\).
Space complexity is also a great burden in the GW calculation,
as manipulating \(\chi_0\) requires to store orbital pair functions
of all occupied and unoccupied orbitals, which costs \(\order{N_{\Ew{e}}^3}\).
The high complexity limits the practical application of GW calculation.

\subsection{Interpolative separable density fitting}
\label{subsec:ISDF}

We use the \emph{Interpolative Separable Density Fitting} (ISDF)
algorithm~\cite{hu2017interpolative} to reduce the time complexity
of evaluating \(\chi_0\).
The ISDF algorithm represents the orbital pair functions \(M_{i,j}(\vec{r})\)
as
\begin{equation}
\label{eq:ISDF}
M_{i,j}(\vec{r})
= \conj{\psi_i(\vec{r})} {\phi_j(\vec{r})}
\thickapprox \sum\limits_{\mu =1}^{N_\mu}
\conj{\psi_i({\vec{r}_{\mu}})} {\phi_j({\vec{r}_\mu})} p_{\mu}(\vec{r}),
\end{equation}
where \(N_{\mu}\) in~\eqref{eq:ISDF} refers to the number of auxiliary
functions.
There are several choices on the interpolation points
\(\{\vec{r}_{\mu}\}_{\mu = 1}^{N_{\mu}}\)
and auxiliary functions \(\{p_{\mu}(\vec{r})\}_{\mu = 1}^{N_{\mu}}\) in the
ISDF algorithm, such as QRCP and
\(k\)-means~\cite{Dong2018ISDFHybridFunctional, lu2015compression}.
We represent~\eqref{eq:ISDF} in matrix form as
\begin{equation}
\label{eq:MatrixExpressionISDF}
M = PC,
\end{equation}
where
\[
M_{(n; i,j)} = \conj{\psi_i(\vec{r}_{n})} {\phi_j(\vec{r}_{n})},
\qquad
P_{(n; \mu)} = p_{\mu}(\vec{r}_{n}),
\qquad
C_{(\mu; i,j)} = \conj{\psi_i(\vec{r}_{\mu})} {\phi_j(\vec{r}_{\mu})}.
\]
Numerical experiments suggest that as the number of wavefunctions increases,
\(N_{\mu}\) needs to grow linearly with respect to the number of
wavefunctions to keep the relative error constant~\cite{lu2015compression}.
Base on this observation, we set
\[
N_{\mu} = k_{\mu}\sqrt{N_1N_2},
\]
where \(N_1\) and \(N_2\) are the numbers of wavefunctions in two sets,
respectively, and \(k_{\mu}\) is a small constant of \(\Theta(10^1)\).

Ma et al.\ implemented the ISDF algorithm in GW
calculation~\cite{ma2021GWISDF}.
Their approach is summarized in Algorithm~\ref{algo:GWMethodClassical}.
By using the ISDF algorithm, Algorithm~\ref{algo:GWMethodClassical} greatly
reduces the dominant quartic term in the time complexity.
The memory requirement is also largely reduced.
However, this algorithm requires an explicit expression of the ISDF
coefficient matrix \(C_{\Ew{vc}}\), which consumes \(\order{N_{\Ew{e}}^3}\)
storage.
Moreover, the cubic terms in time complexity still have ultra-high
prefactors since their approach inverts \(\epsilon\) explicitly.

\begin{algorithm}
\caption{Implement the ISDF algorithm on G\(_0\)W\(_0\) calculation,
from Ma et al.~\cite{ma2021GWISDF}.}
\label{algo:GWMethodClassical}
\spacing{0.8}
\begin{algorithmic}[1]
\REQUIRE KS--DFT wavefunctions \({\{\psi_i\}}_{i = 1}^{N_{\Ew{v}} + N_{\Ew{c}}}\),
KS--DFT band energies
\(\{\varepsilon_i^{\Ew{KS}}\}_{i=1}^{N_{\Ew{v}} + N_{\Ew{c}}}\),
exchange and correlation energies
\(\{{(V_{\Ew{xc}})}_i\}_{i=1}^{N_{\Ew{v}} + N_{\Ew{c}}}\),
indices of interpolation points and auxiliary function matrices \(\{P\}\)
from the ISDF algorithm.
\ENSURE GW energy ${\{\varepsilon^{\GW}_i\}}_{i = 1}^{N_{\Ew{v}} +
N_{\Ew{c}}}$.

\STATE Calculate \hfill \(\order{k_{\mu}N_{\Ew{e}}N_{\Ew{r}}^2+k_{\mu}^2N_{\Ew{e}}^4}\)
\[
    \chi_0 = 2M_{\Ew{vc}}\Omega^{-1}M_{\Ew{vc}}\herm,\qquad
    \epsilon = I- V\chi_0.
\]
\STATE Calculate \(\epsilon^{-1}\) by LU decomposition. \hfill \(\order{N_{\Ew{r}}^3}\)
\STATE Calculate screened Coulomb interaction matrix \hfill \(\order{N_{\Ew{r}}^2}\)
\[
    W_{\Ew{V}} = W - V = (\epsilon^{-1} - I)V.
\]
\STATE Calculate\hfill \(\order{k_{\mu}N_{\Ew{e}}N_{\Ew{r}}^2}\)
\[
    W_{\Ew{vn}} = P_{\Ew{vn}}\herm W_{\Ew{V}} P_{\Ew{vn}},\qquad
    W_{\Ew{nn}} = P_{\Ew{nn}}\herm W_{\Ew{V}} P_{\Ew{nn}},\qquad
    V_{\Ew{vn}} = P_{\Ew{vn}}\herm V    P_{\Ew{vn}}.
\]
\STATE Calculate the part of self-energies induced from the screened exchange
interaction
\hfill \(\order{k_{\mu}^2N_{\Ew{e}}^3}\)
\begin{align*}
\Sigma_{\Ew{SEX\_X}}
={}& \sum_{i\in\occ}
(\conj{\psi_i^{\Ew{vn}}})\herm
\Bigl( W_{\Ew{vn}} \odot \bigl(\Psi^{\Ew{vn}}{(\Psi^{\Ew{vn}})}\herm\bigr)\Bigr)
(\conj{\psi_i^{\Ew{vn}}}),\\
\Sigma_{\Ew{X}}
={}& \sum_{i\in\occ}
(\conj{\psi_i^{\Ew{vn}}})\herm
\Bigl( V_{\Ew{vn}} \odot \bigl(\Psi^{\Ew{vn}}{(\Psi^{\Ew{vn}})}\herm\bigr)\Bigr)
(\conj{\psi_i^{\Ew{vn}}}).
\end{align*}
\STATE Calculate the part of self-energies induced from the Coulomb hole
interaction
\hfill \(\order{k_{\mu}^2N_{\Ew{e}}^3}\)
\begin{align*}
\Sigma_{\Ew{COH}}
={}& \sum_{i\in\all}
(\conj{\psi_i^{\Ew{nn}}})\herm
\Bigl( W_{\Ew{nn}} \odot \bigl(\Psi^{\Ew{nn}}{(\Psi^{\Ew{nn}})}\herm\bigr)\Bigr)
(\conj{\psi_i^{\Ew{nn}}}).
\end{align*}
\STATE Calculate self-energies
\begin{align*}
\Sigma =
\Re\bigl(\diag({\Sigma_{\Ew{SEX\_X}}+\Sigma_{\Ew{X}}+\Sigma_{\Ew{COH}}})\bigr).
\end{align*}
\STATE Calculate quasiparticle energies under G\(_0\)W\(_0\) approximation
\[
    \varepsilon^{\GW}_i =
    \varepsilon_i^{\Ew{KS}} + \Sigma_i - {(V_{\Ew{xc}})}_i.
\]
\end{algorithmic}
\end{algorithm}

\subsection{Sherman--Morrison--Woodbury (SMW) formula combined with ISDF}
\label{subsec:SMWformula}
Suppose that both \(A\) and \(A+UV\herm\) are \(n\times n\) nonsingular
matrices, where \(U\), \(V \in \C^{n\times k}\).
Then the \emph{Sherman--Morrison--Woodbury (SMW)
formula}~\cite{Hager1989UpdateInverseMat} reads
\begin{equation}
   (A+UV\herm)^{-1} = A^{-1} - A^{-1}U(I + V\herm A^{-1}U)^{-1}V\herm A^{-1}.
\label{eq:S--M--W}
\end{equation}

Liu et al.\ in~\cite{liu2015numerical}
suggest that by SMW formula, the inverse of dielectric matrix,
\(\epsilon^{-1}(\omega)\), can be decomposed into a sum of rank-\(1\) matrices
as follows
\begin{equation}
\epsilon^{-1}(\omega)
= I + \sum_{l=1}^{N_{\Ew{v}}N_{\Ew{c}}}
\frac{1}{2\sqrt{\lambda_l}}V(M_{\Ew{vc}}s_l)
\Bigl(\frac{1}{\omega - \sqrt{\lambda_l}+\mi\eta} -
\frac{1}{\omega+ \sqrt{\lambda_l} + \mi\eta}\Bigr)(M_{\Ew{vc}}s_l)\herm,
\label{eq:Liu_dielectric}
\end{equation}
where \(s_l=\Omega^{1/2}u_l\), \(\Omega\) is the same as in~\eqref{eq:Omega},
and \(\{(\lambda_l, u_l)\}_{l=1}^{N_{\Ew{v}}N_{\Ew{c}}}\) are the eigenpairs
of \(\Omega^2+\Omega^{1/2}M\herm VM\Omega^{1/2}\).

Once the eigenpairs of the system described above are known, we can write the
screened exchange part and the Coulomb hole part of the self-energy matrices
as
\[
\begin{aligned}
\Sigma_{\Ew{SEX\_X}}= -\sum_{i\in\occ} (\psi_i \psi_i\herm )
\odot \Biggl(\sum_{l=1}^{N_{\Ew{v}} N_{\Ew{c}}}
\frac{1}{2 \sqrt{\lambda_l}}
\frac{(V M_{\Ew{vc}} s_{l})(V M_{\Ew{vc}} s_{l})\herm}
{\omega - \varepsilon_{i}^{\Ew{KS}} - \sqrt{\lambda_l}}
\Biggr),\\
\Sigma_{\Ew{COH}}= \frac{1}{2}\sum_{i\in\all} (\psi_i \psi_i\herm)
\odot \Biggl(\sum_{l=1}^{N_{\Ew{v}} N_{\Ew{c}}}
\frac{1}{2 \sqrt{\lambda_l}}
\frac{(V M_{\Ew{vc}} s_{l})(V M_{\Ew{vc}} s_{l})\herm}
{\omega - \varepsilon_{i}^{\Ew{KS}} - \sqrt{\lambda_l}}
\Biggr).
\end{aligned}
\]
However, \eqref{eq:Liu_dielectric} requires all eigenvectors of an
\(N_{\Ew{v}}N_{\Ew{c}} \times N_{\Ew{v}}N_{\Ew{c}}\) dense Hermite matrix,
resulting in a time complexity of \(\order{N_{\Ew{e}}^6}\).
This is too expensive for practical computation.

In order to develop a cubic scaling algorithm for computing
\(\epsilon^{-1}\), we propose an alternative approach to apply the SMW
formula.
To simplify the notation, we does not distinguish the operators and their
matrix representations obtained through discretization.
Superscripts on the energies \(\varepsilon\)'s are also omitted.

Note that the ISDF algorithm expresses the orbital pair functions matrix
\(M{_{\Ew{vc}}}\) as~\eqref{eq:MatrixExpressionISDF}.
We substitute~\eqref{eq:MatrixExpressionISDF} into the expression of the 
polarizability matrix \(\chi_0\) and the dielectric matrix \(\epsilon\)
in~\eqref{eq:matrices}, and obtain the matrix representations of \(\chi_0\)
and \(\epsilon\) as
\[\begin{aligned}
\chi_0 &= 2P{_{\Ew{vc}}}C{_{\Ew{vc}}} \Omega^{-1}
(P{_{\Ew{vc}}}C{_{\Ew{vc}}})\herm,\\
\epsilon &= I - 2V(P{_{\Ew{vc}}}C{_{\Ew{vc}}})\Omega^{-1}
(P{_{\Ew{vc}}}C{_{\Ew{vc}}})\herm.
\end{aligned}\]
Applying the SMW formula to invert the dielectric matrix \(\epsilon\), we have
\begin{equation}
\epsilon^{-1} = I + VP_{\Ew{vc}}\Bigl(\frac{1}{2}(\coc)^{-1}
- P_{\Ew{vc}}\herm VP_{\Ew{vc}}\Bigr)^{-1}P_{\Ew{vc}}\herm.
\label{eq:epsilon_with_SMW}
\end{equation}
Clearly~\eqref{eq:epsilon_with_SMW} is much cheaper than the direct inversion
process, because only an \(N{_{\Ew{vc}}}\times N{_{\Ew{vc}}}\) matrix needs to
be explicitly inverted.
Additionally, explicit formulation of \(\chi_0\) and \(\epsilon\) is avoided.
Therefore, \eqref{eq:epsilon_with_SMW} significantly reduces the time
complexity of the inversion process.

We remark that directly applying SMW as~\eqref{eq:epsilon_with_SMW} is still
not fully satisfactory.
The ISDF algorithm provides low-rank properties to the screened Coulomb
interaction operator \(W\) and the self-energy operator \(\Sigma\).
The low-rank properties are not exploited in the conventional G\(_0\)W\(_0\)
calculation strategy.
Even if~\eqref{eq:epsilon_with_SMW} is adopted, further computation in the
G\(_0\)W\(_0\) calculation still involves \(N_{\Ew{r}} \times N_{\Ew{r}}\)
matrices.
Also, we need an explicit expression of the ISDF coefficient matrix
\(C_{\Ew{vc}}\) in~\eqref{eq:epsilon_with_SMW}, which is
$N_{\mu} \times N_{\Ew{v}}N_{\Ew{c}}$.
As inverting \(\epsilon\) requires calculating \(\coc\) first,
our calculation still costs the time complexity of \(\order{N_\Ew{e}^4}\)
and space complexity of \(\order{N_\Ew{e}^3}\).
Therefore, we shall tackle these issues in subsequent subsections.

\subsection{Cauchy integral representation}
In~\cite{Lu2017CubicSA}, Lu and Thicke suggest that the time complexity of
computing \(\chi_0\) can be reduced from \(\order{N_{\Ew{e}}^4}\) to
\(\order{N_{\Ew{e}}^3}\) by using Cauchy integral.
Inspired by Lu and Thicke's work, we propose an algorithm to calculate \(\coc\)
with a cubic time complexity, as listed in Algorithm~\ref{algo:COmegaC}.
The derivation of Algorithm~\ref{algo:COmegaC} is explained in detail in
Supplementary information.

\begin{algorithm}[!tb]
\caption{Calculate~\(\coc\).}
\label{algo:COmegaC}
\spacing{1.0}
\begin{algorithmic}[1]
\REQUIRE A matrix \(\Psi\) that consists of wavefunctions,
indices of interpolative points from the ISDF algorithm,
energies \(\{\varepsilon_i\}_{i=1}^{n}\),
a threshold \(\delta_{\Ew{rel}}\).
\ENSURE Coefficient matrix \(\coc\).
\STATE Calculate
\begin{equation*}
r ={} \frac{\sqrt{Q/q} - 1}{\sqrt{Q/q} + 1},
\end{equation*}
where \(q = \varepsilon_{N_{\Ew{v}}+1} - \varepsilon_{N_{\Ew{v}}}\) and
\(Q = \varepsilon_{N_{\Ew{v}}+N_{\Ew{c}}} - \varepsilon_{N_{\Ew{v}}}\).
\STATE Calculate
\[
R = \int_0^1\frac{\md t}{\sqrt{(1-t^2)(1-r^2t^2)}},
\qquad
L = \frac12\int_0^{1/r}\frac{\md t}{\sqrt{(1+t^2)(1+r^2t^2)}}.
\]
\STATE Use an adaptive quadrature rule to calculate
\[
\coc = -\frac{\sqrt{Qq}}{\pi\mi }\int_{-R+\mi L}^{R+\mi L}\mathrm{Im}
\bigl(J(z(t))\bigr)
\frac{r^{-1}\cn(t)\dn(t)}{(r^{-1} - \sn(t))^2}\md t
\]
until the relative error is below the threshold \(\delta_{\Ew{rel}}\),
where \(J\bigl(z\bigr)\) is given in~(23) of Supplementary information, and 
\[
z(t) = \sqrt{qQ}\biggl(\frac{r^{-1} + \sn(t, r)}{r^{-1} - \sn(t, r)}\biggr) +
\epsilon_{N_{\mathrm{v}}}.
\]
\end{algorithmic}
\end{algorithm}

\subsection{Cubic scaling algorithm for G\(_0\)W\(_0\) calculation}
\label{subsec:cubicGW}
We propose a more efficient and concise version of the G\(_0\)W\(_0\)
calculation strategy using our fast inversion strategy.
Let
\[
K={\Bigl(\frac12{(\coc)}^{-1}-P_{\Ew{vc}}\herm VP_{\Ew{vc}}\Bigr)}^{-1}.
\]
Then~\eqref{eq:epsilon_with_SMW} becomes
\[
\epsilon^{-1} ={} I + VP_{\Ew{vc}} K P_{\Ew{vc}}\herm.
\]
Applying \(\epsilon^{-1}\) in~\eqref{eq:epsilon_with_SMW} to matrices
\(W{_{\Ew{vc}}}\), \(W_{\Ew{vn}}\), and \(W_{\Ew{nn}}\) in steps~3 and~4 of
Algorithm~\ref{algo:GWMethodClassical} yields
\[
V_{\Ew{vn}} = P_{\Ew{vn}}\herm V P_{\Ew{vn}},
\quad
W_{\Ew{vn}} = \bigl(P_{\Ew{vn}}\herm VP_{\Ew{vc}}\bigr) K
\bigl(P_{\Ew{vc}}\herm VP_{\Ew{vn}}\bigr),
\quad
W_{\Ew{nn}} = \bigl(P_{\Ew{nn}}\herm VP_{\Ew{vc}}\bigr) K
\bigl(P_{\Ew{vc}}\herm VP_{\Ew{nn}}\bigr).
\]
Energy matrices under G\(_0\)W\(_0\) approximation in steps~5 and~6 of
Algorithm~\ref{algo:GWMethodClassical} can be viewed as
\begin{equation}
\label{eq:classical56}
\begin{aligned}
\Sigma_{\Ew{SEX\_X}}
&= -\sum_{i\in\occ}(\conj{\psi_i^{\Ew{vn}}}\odot\Psi^{\Ew{vn}})\herm
W_{\Ew{vn}}
(\conj{\psi_i^{\Ew{vn}}}\odot\Psi^{\Ew{vn}}),\\
\Sigma_{\Ew{X}}
&= -\sum_{i\in\occ}
(\conj{\psi_i^{\Ew{vn}}}\odot\Psi^{\Ew{vn}})\herm
V_{\Ew{vn}}
(\conj{\psi_i^{\Ew{vn}}}\odot\Psi^{\Ew{vn}}),\\
\Sigma_{\Ew{COH}}
&= \frac{1}{2}\sum_{i\in\all}
(\conj{\psi_{i}^{\Ew{nn}}}\odot\Psi^{\Ew{nn}})\herm
W_{\Ew{nn}}
(\conj{\psi_{i}^{\Ew{nn}}}\odot\Psi^{\Ew{nn}}).\\
\end{aligned}
\end{equation}
When~\eqref{eq:classical56} is directly adopted for computation, the
complexity is quartic.
The cost can be reduced with the help of the matrix identity
\begin{equation}
\label{eq:Hadamard}
\bigl([\conj a,\conj a,\dotsc,\conj a]\odot B\bigr)\herm C
\bigl([\conj a,\conj a,\dotsc,\conj a]\odot B\bigr)
=B\herm\bigl(C\odot(aa\herm)\bigr)B.
\end{equation}
Then~\eqref{eq:classical56} simplifies to
\[
\begin{aligned}
\Sigma_{\Ew{SEX\_X}}
&= -({\Psi^{\Ew{vn}}})\herm
\Bigl( W_{\Ew{vn}} \odot \bigl(\Psi^{\Ew{vn}}{(\Psi^{\Ew{vn}})}\herm\bigr)\Bigr)
{\Psi^{\Ew{vn}}},\\
\Sigma_{\Ew{X}}
&=-({\Psi^{\Ew{vn}}})\herm
\Bigl( V_{\Ew{vn}} \odot \bigl(\Psi^{\Ew{vn}}{(\Psi^{\Ew{vn}})}\herm\bigr)\Bigr)
\Psi^{\Ew{vn}},\\
\Sigma_{\Ew{COH}}
&=\frac12
(\Psi^{\Ew{nn}})\herm
\Bigl( W_{\Ew{nn}} \odot \bigl(\Psi^{\Ew{nn}}{(\Psi^{\Ew{nn}})}\herm\bigr)\Bigr)
\Psi^{\Ew{nn}},
\end{aligned}
\]
which can be evaluated with a cubic complexity.

In summary, our fast inversion strategy reduces the size of matrices that need
to be inverted from \(N_{\Ew{r}}\times N_{\Ew{r}}\) to \(N{_{\Ew{vc}}}\times
N{_{\Ew{vc}}}\), and decreases the scaling of the matrices to save from a cubic
term plus a quadratic term with a large prefactor to a quadratic term with a
mild prefactor.
We apply our fast inversion strategy to the G\(_0\)W\(_0\) calculation and
rescheduling the multiplication order with~\eqref{eq:Hadamard} to further
exploit the low-rank property.
We present our improved G\(_0\)W\(_0\) calculation strategy in
Algorithm~\ref{algo:GW_method_improved}.
In Algorithm~\ref{algo:GW_method_improved}, the calculation required for
\(W_{\Ew{vn}}\) and \(W_{\Ew{nn}}\) is reduced compared to conventional
G\(_0\)W\(_0\) calculation strategy.
We also mark the time complexity of each step in the right side of the
algorithm and summarize the complexity of all strategies in
Table~\ref{table:TimeSpacecomplexity}.

\begin{algorithm}
\caption{Improved G\(_0\)W\(_0\) calculation based on the ISDF algorithm.}
\label{algo:GW_method_improved}
\small
\spacing{0.8}
\begin{algorithmic}[1]
\REQUIRE KS--DFT wavefunctions \({\{\psi_i\}}_{i = 1}^{N_{\Ew{v}} + N_{\Ew{c}}}\),
KS--DFT band energies
\(\{\varepsilon_i^{\Ew{KS}}\}_{i=1}^{N_{\Ew{v}} + N_{\Ew{c}}}\),
exchange and correlation energies
\(\{{(V_{\Ew{xc}})}_i\}_{i=1}^{N_{\Ew{v}} + N_{\Ew{c}}}\),
indices of the interpolation points and the auxiliary function
matrices~\(\{P\}\) from the ISDF algorithm.
\ENSURE GW energies
\({\{\varepsilon^{\GW}_i\}}_{i = 1}^{N_{\Ew{v}} + N_{\Ew{c}}}\).

\STATE Calculate \(\coc\) by Algorithm~\ref{algo:COmegaC} and invert \(\coc\).
\hfill \(\order{N_{\Ew{e}}^3}\)
\STATE Calculate \hfill \(\order{k_{\mu}^2N_{\Ew{e}}^2N_{\Ew{r}}}\)
\[
P_{\Ew{vc}}\herm V P{_{\Ew{vc}}},\qquad
P_{\Ew{vn}}\herm V P{_{\Ew{vc}}},\qquad
P_{\Ew{nn}}\herm V P{_{\Ew{vc}}}.
\]
\STATE Calculate \hfill \(\order{k_{\mu}^3N_{\Ew{e}}^3}\)
\[
K = \Bigl(\frac{1}{2}(\coc)^{-1} - P_{\Ew{vc}}\herm V P{_{\Ew{vc}}}\Bigr)^{-1}.
\]
\STATE Calculate \hfill \(\order{k_{\mu}^2N_{\Ew{e}}^2N_{\Ew{r}}}\)\\
\[
W_{\Ew{vn}}
= \bigl(P_{\Ew{vn}}\herm V P{_{\Ew{vc}}}\bigr)K\bigl(P_{\Ew{vn}}\herm
VP{_{\Ew{vc}}}\bigr)\herm,
\qquad
W_{\Ew{nn}}
= \bigl(P_{\Ew{nn}}\herm VP{_{\Ew{vc}}}\bigr)K\bigl(P_{\Ew{nn}}\herm
VP{_{\Ew{vc}}}\bigr)\herm.
\]
\STATE Calculate SEX interaction
\hfill \(\order{k_{\mu}^2N_{\Ew{e}}^3}\)
\begin{align*}
\Sigma_{\Ew{SEX\_X}}
&= -({\Psi^{\Ew{vn}}})\herm
\Bigl( W_{\Ew{vn}} \odot \bigl(\Psi^{\Ew{vn}}{(\Psi^{\Ew{vn}})}\herm\bigr)\Bigr)
{\Psi^{\Ew{vn}}},\\
\Sigma_{\Ew{X}}
&=-({\Psi^{\Ew{vn}}})\herm
\Bigl( V_{\Ew{vn}} \odot \bigl(\Psi^{\Ew{vn}}{(\Psi^{\Ew{vn}})}\herm\bigr)\Bigr)
\Psi^{\Ew{vn}},
\end{align*}
\STATE Calculate COH correlation
\hfill \(\order{k_{\mu}^2N_{\Ew{e}}^3}\)
\[
\Sigma_{\Ew{COH}}
=\frac12 (\Psi^{\Ew{nn}})\herm
\Bigl( W_{\Ew{nn}} \odot \bigl(\Psi^{\Ew{nn}}{(\Psi^{\Ew{nn}})}\herm\bigr)\Bigr)
\Psi^{\Ew{nn}}.
\]
\STATE Evaluate the self-energies:
\[
\Sigma = \Re\bigl(\diag({\Sigma_{\Ew{SEX\_X}}+\Sigma_{\Ew{X}}+\Sigma_{\Ew{COH}}})\bigr).
\]
\STATE Evaluate quasiparticle energies
\[
    \varepsilon^{\Ew{GW}}_i = \varepsilon_i^{\Ew{KS}} + {\Sigma}_i - {(V_{\Ew{xc}})}_i.
\]
\end{algorithmic}
\end{algorithm}

\subsection{Error analysis of implementing the ISDF algorithm
on G\(_0\)W\(_0\) calculation}
\label{subsec:errorAnalysis}
Algorithm~\ref{algo:GW_method_improved} introduces quadrature error in the
calculation of \(C_{\Ew{vc}}\Omega^{-1} C_{\Ew{vc}}\herm\) by Cauchy integral
and implementing the ISDF algorithm to decompose the matrices of orbital
pairs.
The quadrature error associated with Cauchy integral can be well-bounded.
It is pointed out in \cite[Lemma~3.1]{Lu2017CubicSA} that the upper bound of
quadrature error introduced in the Cauchy integral is
\begin{equation}
\label{eq:rel_err_Cauchy}
\mathcal{O}\biggl(\exp\Bigl(\frac{-\pi^2 N_\lambda}
{2\log (Q/q)+6}\Bigr)\biggr),
\end{equation}
where \(N_{\lambda}\) is the number of quadrature nodes, and \(Q/q\) is the
condition number determined by the bandgap of the system.
As the error decreases geometrically with respect to the number of quadrature
nodes, we can decrease the error with little cost.

For the error introduced by the ISDF algorithm, we provide the upper bound of
relative error of the self-energies in terms of the error introduced to the
orbital pair functions matrix \(M\) as
\begin{equation}\begin{aligned}
\label{eq:rel_err_ISDF}
\abs{\delta\bigl[\Sigma(n)\bigr]}
\leq{}& \norm{\delta[X]}_2\tr\left({A_{(\all, n; \all, n)}}\right) +
\abs{\tr(\delta[A_{(\all, n; \all, n)}])}\notag\\
&+ \frac{1}{2} \lvert\tr(\delta[A])\rvert+o(\Fnorm{\delta[M]})\notag\\
={}& \norm{\delta[X]}_2\tr\left({A_{(\all, n; \all, n)}}\right) + \frac{3}{2}
\lvert\tr(\delta[A_{(\all, n; \all, n)}])\rvert + o(\Fnorm{\delta[M]}).
\end{aligned}\end{equation}
Here, \(A\) refers to the four-center two-electron integral matrix, and \(\all\) in
the subscripts indicates indices corresponding to all states.
Detailed derivation of~\eqref{eq:rel_err_ISDF} is also provided in
Supplementary information.
The stability of the improved GW strategy is then ensured
by~\eqref{eq:rel_err_Cauchy} and~\eqref{eq:rel_err_ISDF}, as long as the ISDF
algorithm provides a relatively good approximation of \(M\).
Such an assumption is usually plausible, and is easily attainable in practice.

\begin{table}[!tb]
\centering
\caption{List of notation.}
\label{table:notation}
\spacing{1.5}
\begin{tabular}{p{0.2\textwidth}p{0.75\textwidth}}
\hline
Notation                              & Description\\
\hline
\(\psi_i(\vec{r})\), \(\phi_i(\vec{r})\)  & \(i\)-th wavefunction\\
\(\Psi\) & wavefunction matrix, \(\Psi_{(i, k)} = \psi_i(\vec{r}_k)\).\\
\(\varepsilon_i^{\Ew{KS}}\)                       & \(i\)-th energy, from
KS--DFT calculation\\
\(\varepsilon_i^{\Ew{GW}}\)                       & \(i\)-th quasiparticle
energy, from GW calculation\\
\(\occ\), \(\unocc\), \(\all\)              & set of indices of all occupied
states, all unoccupied states and all states\\
\(N_{\Ew{v}}\)                          & number of occupied states\\
\(N_{\Ew{c}}\)                          & number of unoccupied states\\
\(N_{\Ew{e}}\)                          & number of electrons\\
\(N_{\Ew{r}}\)                          & number of discrete points in real
space\\
\(N_{\mu}\)                             & number of auxiliary functions in the
ISDF algorithm\\
\(k_{\mu}\)                             & ISDF coefficient\\
\(M_{\Ew{vc}}\), \(M_{\Ew{vn}}\) and \(M_{\Ew{nn}}\)
& orbital pair functions of occupied and unoccupied states,
  occupied states and all states, all states and all states, respectively\\
\({\cdot}_{\Ew{vc}}\), \({\cdot}_{\Ew{vn}}\) and \({\cdot}_{\Ew{nn}}\)
    & subscripts indicate functions or operators that are obtained 
    from ISDF process of \(M_{\Ew{vc}}\),
    \(M_{\Ew{vn}}\), and \(M_{\Ew{nn}}\), respectively\\
\({\cdot}^{\Ew{vc}}\), \({\cdot}^{\Ew{vn}}\) and \({\cdot}^{\Ew{nn}}\)
    & superscripts indicate wavefunctions on interpolation points from ISDF
    process\\
\(A\trans\)                             & transpose of \(A\)\\
\(A\herm\)                              & conjugate transpose of \(A\)\\
\(\conj{A}\)                            & conjugate of \(A\)\\
\(\odot\)                               & Hadamard product\\
\(\otimes\)                             & Kronecker product\\
\(I\)                                   & identity matrix\\
\hline
\end{tabular}
\end{table}

\begin{table}[!tb]
\centering
\caption{Complexity for different implementations of G\(_0\)W\(_0\) calculation.
In our numerical experiments, \(N_{\Ew{r}} \approx 1000N_{\Ew{e}}\).
}
\label{table:TimeSpacecomplexity}
\footnotesize
\begin{tabular}{ccccccc}
\hline
& \multicolumn{2}{c}{{Conventional (BerkeleyGW)}} &
\multicolumn{2}{c}{{Ma et al.'s work~\cite{ma2021GWISDF}}} &
\multicolumn{2}{c}{{This work}} \\
& Time & Space & Time & Space & Time & Space \\
\hline
\multirow{2}{*}{Operators} &
\multicolumn{2}{c}{\(\chi\), \(\epsilon\), \(\epsilon^{-1}\), \(W\)} &
\multicolumn{2}{c}{\(\chi\), \(\epsilon\), \(\epsilon^{-1}\), \(W\)} &
\multicolumn{2}{c}{\(\coc\), \(K\), \(P_{\Ew{vc}}\herm VP_{\Ew{vc}}\)} \\
& \vphantom{\(a^{a^{a^a}}\)}\(\order{N_{\Ew{r}}^3 + N_{\Ew{e}}^2N_{\Ew{r}}^2}\)
& \(4N_{\Ew{r}}^2\) & \(\order{N_{\Ew{r}}^3 + k_{\mu}^2N_{\Ew{e}}^4}\) &
\(4N_{\Ew{r}}^2\) &
\(\order{k_{\mu}N_{\Ew{e}}^3}\) & \(4k_{\mu}^2N_{\Ew{e}}^2\) \\
\hline
\vphantom{\(a^{a^{a^a}}\)}    {\(\Sigma_{\Ew{SEX\_X}}\), \(\Sigma_{\Ew{X}}\)}
    & \(\order{N_{\Ew{e}}^2N_{\Ew{r}}^2}\)
    & \(N_{\Ew{e}}^2\)
    & \(\order{k_{\mu}N_{\Ew{e}}N_{\Ew{r}}^2}\)
    & \(k_{\mu}^2N_{\Ew{e}}^2\)
    & \(\order{k_{\mu}^2N_{\Ew{e}}^2N_{\Ew{r}}}\)
    & \(k_{\mu}^2{N_{\Ew{e}}^2}\)  \\
    \hline
\vphantom{\(a^{a^{a^a}}\)}    {\(\Sigma_{\Ew{COH}}\)}
    & \(\order{N_{\Ew{e}}^2N_{\Ew{r}}^2}\)
    & \(N_{\Ew{e}}^2\)
    & \(\order{k_{\mu}N_{\Ew{e}}N_{\Ew{r}}^2}\)
    & \(k_{\mu}^2N_{\Ew{e}}^2\)
    & \(\order{k_{\mu}^2N_{\Ew{e}}^2N_{\Ew{r}}}\) & \(k_{\mu}^2N_{\Ew{e}}^2\)\\
\hline
\end{tabular}
\end{table}






\begin{addendum}
\item This work is partly supported by the Innovation Program for Quantum Science and Technology (2021ZD0303306), the Strategic Priority Research Program of the Chinese Academy of Sciences (XDB0450101), the National Natural Science Foundation of China (22288201, 22173093, 21688102), by the Anhui Provincial Key Research and Development Program (2022a05020052), the National Key Research and Development Program of China (2016YFA0200604, 2021YFB0300600), and the CAS Project for Young Scientists in Basic Research (YSBR-005). The authors declare no competing financial interest.

\end{addendum}

\newpage

\footnotesize{
\bibliography{ref}
}

\newpage

\end{document}


\maketitle
\tableofcontents

\newpage

\setcounter{equation}{19}

\appendix 

\section{Parameters of test systems}
\label{sec:parameters}
\begin{table}[!h]
\centering
\caption{Parameters of the systems in numerical examples.}
\label{table:parameters}
\begin{tabular}{lcccc}
\hline
Name       &valence electrons \(N_{\Ew{e}}\)    &occupied states \(N_{\Ew{v}}\) 
           &unoccupied states \(N_{\Ew{c}}\)    & Ecut (Ha)  \\ 
\hline
Si\(_8\)     &\(32\)            &\(16\)          &\(16\)         & \(10.0\)   \\
Si\(_{16}\)  &\(64\)            &\(32\)          &\(32\)         & \(10.0\)   \\
Si\(_{24}\)  &\(96\)            &\(48\)          &\(48\)         & \(10.0\)   \\
Si\(_{32}\)  &\(128\)           &\(64\)          &\(64\)         & \(10.0\)   \\
Si\(_{64}\)  &\(256\)           &\(128\)         &\(128\)        & \(10.0\)   \\
(SrTiO\(_3\))\(_8\)&\(192\)       &\(96\)          &\(96\)         & \(20.0\) \\
C\(_{60}\)   &\(240\)           &\(120\)         &\(120\)        & \(10.0\)   \\
\hline
\end{tabular}
\end{table}

\newpage
\section{Detailed derivation}
\label{sec:appendix}
\subsection{Detailed derivation process of Algorithm~2}
\label{subsec:detailedCauchyIntegral}
The time complexity of calculating the coefficient coupling matrix \(\coc\) is
quartic, because the \((\mu,\nu)\) entry of \(\coc\) is of the form
\begin{equation}
\label{eq:COmegaCstar}
(\coc)_{(\mu,\nu)}
=\sum_{\substack{i\in\occ\\j\in\unocc}}
\frac{\conj{\psi_i(\vec{r}_{\mu})} \psi_i(\vec{r}_{\nu})
\psi_j(\vec{r}_{\mu}) \conj{\psi_j(\vec{r}_{\nu})}}
{\varepsilon_i - \varepsilon_j}.
\end{equation}

To simplify~\eqref{eq:COmegaCstar}, we select a closed contour \(\mathfrak{I}\)
in the complex plane, so that all eigenvalues correspond to unoccupied states
are enclosed while none of the eigenvalues correspond to occupied states are
enclosed.
Then, the coefficients in~\eqref{eq:COmegaCstar} can be expressed as
\begin{equation}
\label{eq:coeffchi0}
\frac{1}{\varepsilon_i - \varepsilon_j} = -\frac{1}{2\pi\mi
}\int_{\mathfrak{I}}\frac{1}{(\lambda - \varepsilon_i)(\lambda -
\varepsilon_j)}\md\lambda.
\end{equation}
By making use of~\eqref{eq:coeffchi0}, we obtain
\begin{align}
{(\coc)}_{(\mu,\nu)}
={}& -\frac{1}{2\pi\mi}\sum_{\substack{i\in\occ\\j\in\unocc}}
\conj{\psi_i(\vec{r}_{\mu})}\psi_i(\vec{r}_{\nu})
\psi_j(\vec{r}_{\mu})\conj{\psi_j(\vec{r}_{\nu})}
\int_{\mathfrak{I}}\frac{1}{(\lambda - \varepsilon_i)(\lambda - \varepsilon_j)}
\md\lambda\nonumber\\
={}& -\frac{1}{2\pi\mi}\sum_{\substack{i\in\occ\\j\in\unocc}}
\int_{\mathfrak{I}}\biggl(
\frac{\conj{\psi_i(\vec{r}_{\mu})}\psi_i(\vec{r}_{\nu})}{\lambda - \varepsilon_i}
\biggr)
\biggl(
\frac{\psi_j(\vec{r}_{\mu})\conj{\psi_j(\vec{r}_{\nu})}}
{\lambda - \varepsilon_j}
\biggr)\md\lambda.
\label{eq:SinglePointCOmegaC}
\end{align}

Let us introduce two matrices \(\Psi_{\Ew{v}}\) and \(\Psi_{\Ew{c}}\) that
consist of occupied and unoccupied wavefunctions on interpolation
points, respectively.
These matrices are defined through
\begin{align*}
{(\Psi_{\Ew{v}})}_{({\mu}, i)} = \psi_i(\vec{r}_{\mu}),
\qquad \forall\, i\in\occ, \\
{(\Psi_{\Ew{c}})}_{({\mu}, j)} = \psi_j(\vec{r}_{\mu}),
\qquad \forall\, j \in \unocc.
\end{align*}
We then reformulate~\eqref{eq:SinglePointCOmegaC} in matrix expression as
\begin{align*}
&{\coc} \\
={}& -\frac{1}{2\pi\mi}\sum_{\substack{i\in\occ\\j\in\unocc}}
\int_{\mathfrak{I}}
\Bigl(\conj{{({\Psi_{\Ew{v}}})}_{(:,i)}}
\frac{1}{\lambda - \varepsilon_i}
\bigl({({\Psi_{\Ew{v}}})}_{(:,i)}\bigr)\trans \Bigr)
\Bigl({({\Psi_{\Ew{c}}})}_{(:,j)}
\frac{1}{\lambda - \varepsilon_j}
\bigl(\conj{{(\Psi_{\Ew{c}})}_{(:,j)}}\bigr)\trans\Bigr)
\md\lambda\\
={}& -\frac{1}{2\pi\mi}\int_{\mathfrak{I}}
\Bigl(\sum_{i\in\occ}\conj{{({\Psi_{\Ew{v}}})}_{(:,i)}}
\frac{1}{\lambda - \varepsilon_i}
\bigl({({\Psi_{\Ew{v}}})}_{(:,i)}\bigr)\trans \Bigr)
\Bigl(\sum_{j\in\unocc}{({\Psi_{\Ew{c}}})}_{(:,j)}
\frac{1}{\lambda - \varepsilon_j}
\bigl(\conj{{(\Psi_{\Ew{c}})}_{(:,j)}}\bigr)\trans\Bigr)
\md\lambda\\
={}&-\frac{1}{2\pi\mi}\int_{\mathfrak{I}}
\Bigl(\conj{\Psi_{\Ew{v}}}
\Omega_{\Ew{v}}^{-1}(\lambda)
(\conj{\Psi_{\Ew{v}}})\herm\Bigr)
\odot \Bigl({\Psi_{\Ew{c}}} \Omega_c^{-1}(\lambda)
({\Psi_{\Ew{c}}})\herm \Bigr)
\md\lambda\\
={}& -\frac{1}{2\pi\mi}\int_{\mathfrak{I}}J(\lambda)\md\lambda,
\end{align*}
where \({\Omega_{\Ew{v}}}\) and \({\Omega_{\Ew{c}}}\) are diagonal matrices
defined as
\Be{
{\Omega_{\Ew{v}}}_{(i, i)} = {\lambda - \varepsilon_i},
    \qquad \forall\, i\in\occ,\\
{\Omega_{\Ew{c}}}_{(j, j)} = {\lambda - \varepsilon_j},
    \qquad \forall\, j\in\unocc,
}
and
\begin{equation}
\label{eq:defJlambda}
J(\lambda) ={} \Bigl(\conj{\Psi_{\Ew{v}}}\Omega_{\Ew{v}}^{-1}(\lambda)
(\conj{\Psi_{\Ew{v}}})\herm\Bigr) \odot
\Bigl({\Psi_{\Ew{c}}}\Omega_{\Ew{c}}^{-1}(\lambda)({\Psi_{\Ew{c}}})\herm \Bigr).
\end{equation}
It is worth noting that the time reversal symmetric
property~(8) implies
\begin{align*}
{\bigl(\conj{J(\conj{\lambda})}\bigr)}_{(\mu,\nu)}
={}& \sum_{\substack{i\in\occ\\j\in\unocc}}
\conj{\Bigl(\conj{{\Psi_{\Ew{v}}}_{(\mu,i)}}
\frac{1}{\conj{\lambda} - \varepsilon_i}
{\Psi_{\Ew{v}}}_{(\nu,i)}\Bigr)}
\conj{\Bigl({\Psi_{\Ew{c}}}_{(\mu,j)}
\frac{1}{\conj{\lambda}-\varepsilon_j}
\conj{{\Psi_{\Ew{c}}}_{(\nu,j)}}\Bigr)}\\
={}& \sum_{\substack{i\in\occ\\j\in\unocc}} 
\Bigl({\Psi_{\Ew{v}}}_{(\mu,i)}
\frac{1}{\lambda - \varepsilon_i}
\conj{\Psi_{\Ew{v}}}_{(\nu,i)} \Bigr)
\Bigl(\conj{{\Psi_{\Ew{c}}}_{(\mu,j)}}
\frac{1}{{\lambda} - \varepsilon_j}
{\Psi_{\Ew{c}}}_{(\nu,j)}\Bigr)\\
={}& \sum_{\substack{i\in\occ\\j\in\unocc}}
\Bigl(\conj{{\Psi_{\Ew{v}}}_{(\mu,i)}}
\frac{1}{\lambda-\varepsilon_i}
{\Psi_{\Ew{v}}}_{(\nu,i)} \Bigr)
\Bigl({\Psi_{\Ew{c}}}_{(\mu,j)}
\frac{1}{{\lambda} - \varepsilon_j}
{\conj{\Psi_{\Ew{c}}}_{(\nu,j)}}\Bigr)\\
={}& {\bigl(J(\lambda)\bigr)}_{(\mu,\nu)}.
\end{align*}

By choosing the contour \(\mathfrak{I}\) to be symmetric with respect to the
real axis, we obtain
\[
\coc=-\frac{1}{2\pi\mi}\int_{\mathfrak{I}}J(\lambda)\md\lambda
=-\frac{1}{\pi\mi}\int_{\mathfrak{I}_{+}}\Im\bigl(J(\lambda)\bigr)
\md\lambda,
\]
where \(\mathfrak{I}_{+}\) is the segment of \(\mathfrak{I}\) lies in the
upper half plane.
It remains to choose an appropriate contour.

Lu and Thicke in~\cite{Lu2017CubicSA} suggest that elliptic integral of the
first kind is a good choice of integration path with well robustness.
Here, we provide a brief introduction to this approach.
Define \(q = \varepsilon_{N_{\Ew{v}}+1} - \varepsilon_{N_{\Ew{v}}}\),
\(Q = \varepsilon_{N_{\Ew{v}}+N_{\Ew{c}}} - \varepsilon_{N_{\Ew{v}}}\),
and
\[
r = \frac{\sqrt{Q/q} - 1}{\sqrt{Q/q}+1},
\qquad
R = \int_0^1\frac{\md t}{\sqrt{(1-t^2)(1-r^2t^2)}},
\qquad
L = \frac{1}{2}\int_0^{1/r}\frac{\md t}{\sqrt{(1+t^2)(1+r^2t^2)}}.
\]
We map the line segment with vertices \(-R+\mi L\) and \(R+\mi L\) to
\(\mathfrak{I}_+\) via
\[
z = \sqrt{qQ}\biggl(\frac{r^{-1} + \sn(t, r)}{r^{-1} - \sn(t, r)}\biggr) +
\epsilon_{N_{\mathrm{v}}},
\]
where \(\sn(t,r)\) refers to the Jacobi elliptic function.

\subsection{Detailed error analysis}
\label{subsec:Detailederroranalysis}
Let \(M = M_{\Ew{vc}}\) in the following paragraphs of this section to
simplify notation.
We also use \(\occ\), \(\unocc\), and \(\all\) as subscripts.
For example, \(M_{\occ\all}\) refers to the orbital pair functions matrix of
occupied and all states, which is equal to \(M_{\Ew{vn}}\) in the sections above.

To derivate the error bound given in~(19), we first
consider the conventional G\(_0\)W\(_0\) calculation strategy.
Combining~(5) and~(7) yields
\begin{equation}
\label{eq:Sigma_G0W0_2}
\begin{aligned}
\Sigma_{\Ew{X}}(n)
={}& -\sum_{i\in\occ}\sum_{\vec{r}, \vec{r}^{\prime}}
\psi_{n}(\vec{r})\conj{\psi_i(\vec{r})}
V(\vec{r},\vec{r}^{\prime})
\conj{\psi_{n}(\vec{r}^{\prime})}\psi_i(\vec{r}^{\prime}),\\
\Sigma_{\Ew{SEX\_X}}(n)
={}& -\sum_{i\in\occ}\sum_{\vec{r}, \vec{r}^{\prime}}
\psi_{n}(\vec{r})\conj{\psi_i(\vec{r})}
\bigl(W(\vec{r},\vec{r}^{\prime}) - V(\vec{r}, \vec{r}^{\prime})\bigr)
\conj{\psi_{n}(\vec{r}^{\prime})}\psi_i(\vec{r}^{\prime}),\\
\Sigma_{\Ew{COH}}(n)
={}& \frac{1}{2}\sum_{i\in \all}\sum_{\vec{r}, \vec{r}^{\prime}}
\psi_{n}(\vec{r})\conj{\psi_{i}(\vec{r})}
\bigl(W(\vec{r},\vec{r}^{\prime}) - V(\vec{r}, \vec{r}^{\prime})\bigr)
\conj{\psi_{n}(\vec{r}^{\prime})}\psi_{i}(\vec{r}^{\prime}).
\end{aligned}
\end{equation}
According to the SMW formula, we have
\[
\epsilon^{-1} = I + VM\Bigl(\frac{1}{2}\Omega - M\herm V M\Bigr)^{-1} M\herm.
\]
Substituting \(W\) in~\eqref{eq:Sigma_G0W0_2} with~(4),
we rewrite~\eqref{eq:Sigma_G0W0_2} as
\begin{equation}
\begin{aligned}
\label{eq:Sigma_G0W0_3}
\Sigma_{\Ew{SEX\_X}}(n)
={}& -\sum_{i\in\occ}\sum_{\vec{r}, \vec{r}^{\prime}}
\conj{\psi_{n}(\vec{r})}{\psi_i(\vec{r})}
\Bigl(V M\Bigl(\frac{1}{2}\Omega - M\herm V M\Bigr)^{-1}M\herm V \Bigr)
(\vec{r},\vec{r}^{\prime})
{\psi_{n}(\vec{r}^{\prime})}\conj{\psi_i(\vec{r}^{\prime})},\\
\Sigma_{\Ew{COH}}(n) ={}& \frac{1}{2} \sum_{i\in \all}\sum_{\vec{r}, \vec{r}^{\prime}}
\conj{\psi_{n}(\vec{r})}{\psi_{i}(\vec{r})}
\Bigl(V M\Bigl(\frac{1}{2}\Omega - M\herm V M\Bigr)^{-1}M\herm V\Bigr)
(\vec{r}, \vec{r}^{\prime})
{\psi_{n}(\vec{r})}\conj{\psi_{i}(\vec{r})},\\
\Sigma_{\Ew{X}}(n)
={}& -\sum_{i\in\occ}\sum_{\vec{r}, \vec{r}^{\prime}}
\conj{\psi_{n}(\vec{r})}{\psi_i(\vec{r})}
V(\vec{r},\vec{r}^{\prime})
{\psi_{n}(\vec{r}^{\prime})}\conj{\psi_i(\vec{r}^{\prime})}.\\
\end{aligned}
\end{equation}
Let \(A\) be the four-center two-electron repulsion integral matrix of all orbitals,
that is
\[\begin{aligned}
A_{(ij, kl)} ={}&
\int_{\R^3 \times \R^3}
\psi_{i}(\vec{r})\conj{\psi_{j}(\vec{r})}
V(\vec{r},\vec{r}^{\prime})
\conj{\psi_{k}(\vec{r}^{\prime})}\psi_{l}(\vec{r}')\md \vec{r} \md\vec{r}^{\prime}
= \conj{A_{(kl, ij)}},\qquad
i,j,k,l \in \all.
\end{aligned}\]
The last equality is guaranteed by the symmetry of \(V(\vec{r}, \vec{r})\).
By representing \(A\) with the orbital pair functions matrix \(M\) as
\[
A = M\herm V M,
\]
equation~\eqref{eq:Sigma_G0W0_3} is reformulated as
\begin{equation}\begin{aligned}
\label{eq:Sigma_by_I4c2e}
\Sigma_{\Ew{SEX\_X}}(n) ={}& -\sum_{i\in\occ}
A_{(\occ\unocc, in)}\herm
\Bigl(\frac{1}{2}\Omega - A_{(\occ\unocc,\occ\unocc)}\Bigr)^{-1}
A_{(\occ\unocc,in)},\\
\Sigma_{\Ew{COH}}(n) ={}& \frac{1}{2}\sum_{n'\in\all}
A_{(\occ\unocc,n'n)}\herm
\Bigl(\frac{1}{2}\Omega - A_{(\occ\unocc, \occ\unocc)}\Bigr)^{-1}
A_{(\occ\unocc,n'n)},\\
\Sigma_{\Ew{X}}(n) ={}& -\sum_{i\in\occ}
A_{(in, in)}.
\end{aligned}\end{equation}
Summing these terms together, the self-energies become
\begin{equation}
\begin{aligned}
\Sigma(n) ={}&\Sigma_{\Ew{SEX}}(n)+\Sigma_{\Ew{COH}}(n)+\Sigma_{\Ew{X}}(n)\\
={}& \frac{1}{2}\sum_{j\in\unocc}A_{(\occ\unocc, nj)}\herm
\Bigl(\frac{1}{2}\Omega - A_{(\occ\unocc, \occ\unocc)}\Bigr)^{-1}
A_{(\occ\unocc, jn)}\\
&-
\frac{1}{2}\sum_{i\in\occ}
A_{(\occ\unocc,in)}\herm
\Bigl(\frac{1}{2}\Omega - A_{(\occ\unocc,\occ\unocc)}\Bigr)^{-1}
A_{(\occ\unocc,in)}
- \sum_{i\in\occ} A_{(in,in)},
\label{eq:SigmawithA}
\end{aligned}
\end{equation}
which are solely expressed in terms of the matrix \(A\).

Let us introduce a positive semi-definite matrix
\(X = V^{1/2} M (-2\Omega^{-1}) M\herm V^{1/2}\).
The summation in~\eqref{eq:SigmawithA} can be regarded as a summation of
quadratic forms with the same coefficient matrix
\(V^{1/2} M \bigl(\frac{1}{2}\Omega - M\herm V M\bigr)^{-1}M\herm V^{1/2}\).
Let us rewrite the coefficient matrix as
\begin{align}
& V^{1/2} M\Bigl(\frac{1}{2}\Omega - M\herm V M\Bigr)^{-1}M\herm V^{1/2}
\nonumber \\
={}& \bigl(V^{1/2} M\bigr)
\Bigl(2\Omega^{-1} + 4\Omega^{-1}\bigl(V^{1/2} M\bigr)\herm
\Bigl(I - 2\bigl(V^{1/2} M\bigr)\Omega^{-1}
\bigl(V^{1/2} M\bigr)\herm)\Bigr)^{-1}
\bigl(V^{1/2} M\bigr)\Omega^{-1}\Bigr)
\bigl(V^{1/2} M\bigr)\herm\nonumber\\
={}& -X + X(I + X)^{-1}X\nonumber\\
={}& -X(I + X)^{-1}\nonumber\\
={}& 2V^{1/2} M \Omega^{-1} M\herm V^{1/2}(I - 2V^{1/2}
M \Omega^{-1} M\herm V^{1/2})^{-1}\label{eq:KernelSigma}.
\end{align}
We take the \(\Sigma_{\Ew{SEX}}\) in~\eqref{eq:Sigma_by_I4c2e} as an example;
\(\Sigma_{\Ew{COH}}\) is computed in the same manner.
Using the definition of \(A\) and~\eqref{eq:KernelSigma}, we reformulate 
\(\Sigma_{\Ew{SEX}}\) with trace notation as
\begin{align*}
\Sigma_{\Ew{SEX}}(n)
&= \sum_{i\in\occ}
A_{(\occ\unocc, in)}\herm
\Bigl(\frac12\Omega - A_{(\occ\unocc, \occ\unocc)}\Bigr)^{-1}
A_{(\occ \unocc, in)}\notag\\
&=\tr\left(M_{\occ n}\herm V^{1/2}\Bigl(V^{1/2} M
\bigl(\frac{1}{2}\Omega - M\herm V M\bigr)^{-1}
M\herm V^{1/2}\Bigr)V^{1/2}M_{\occ n}\right)\notag\\
&=\tr\left(2M_{\occ n}\herm V^{1/2}V^{1/2} M \Omega^{-1} M\herm V^{1/2}
\left(I - 2V^{1/2} M \Omega^{-1} M\herm V^{1/2}\right)^{-1}
V^{1/2}M_{\occ n}\right)\notag\\
&= \tr\left(\left(2V^{1/2} M \Omega^{-1} M\herm
V^{1/2} (I - 2V^{1/2}M\Omega^{-1}M\herm V^{1/2})^{-1}\right)
V^{1/2}M_{\occ n}M_{\occ n}\herm V^{1/2}\right)\notag\\
&= -\tr\left(X(I+X)^{-1}V^{1/2}M_{\occ n}M_{\occ n}\herm V^{1/2}\right).
\end{align*}
In the same manner, the self-energies in~\eqref{eq:SigmawithA} become
\begin{align*}
\Sigma(n)
={}&-\frac{1}{2}\tr(X(I + X)^{-1}V^{1/2} M_{\unocc n}M_{\unocc n}\herm V^{1/2})
- \tr(A_{(\occ n, \occ n)})\notag\\
{}&+\frac{1}{2}\tr(X(I + X)^{-1}V^{1/2} M_{\occ n}M_{\occ n}\herm V^{1/2}).
\end{align*}
Differentiating the self-energies with respect to \(M\) leads to
\begin{align}
\abs{\delta\bigl[\Sigma(n)\bigr]}
\leq{}& \frac{1}{2} \abs{\tr\Bigl(X(I + X)^{-1} \delta\bigl[V^{1/2} M_{\all
n}M_{\all n}\herm V^{1/2}\bigr] \Bigr)}
+ \abs{\tr\bigl(\delta [A]_{(\occ n, \occ n)}\bigr)}\nonumber\\
&+ \frac{1}{2}\abs{\tr\Bigl(\delta\bigl[X(I + X)^{-1}\bigr]
V^{1/2} M_{\all n}M_{\all n}\herm V^{1/2} \Bigr)} + o(\Fnorm{\delta[M]}).
\label{eq:EstimationDeltaSigma1}
\end{align}
Notice that
\begin{align*}
\tr(A_{(\all n,\all n)})
&= \tr\left(V^{1/2} M_{\all n}M_{\all n}\herm V^{1/2}\right),\\
\norm{X (I+X)^{-1}}_2 &\leq 1, \qquad\norm{(I+X)^{-1}}_2\leq 1,\\
\norm{\delta\bigl[X (I+X)^{-1}\bigr]}_2
&\leq{} \norm{\delta[X](I + X)^{-1}}_2 + \norm{X(I+X)^{-1}\delta[X](I+X)^{-1}}_2\\
&\leq \norm{\delta[X]}_2\norm{(I + X)^{-1}}_2
+ \norm{X(I+X)^{-1}}_2\norm{(I+X)^{-1}}_2\norm{\delta[X]}_2\\
&\leq 2\norm{\delta[X]}_2.
\end{align*}
Substituting these inequalities to~\eqref{eq:EstimationDeltaSigma1}, we
conclude that
\begin{align*}
\abs{\delta\bigl[\Sigma(n)\bigr]}
\leq{}&\frac{1}{2} \abs{\tr\Bigl(\delta\bigl[V^{1/2}M_{\Ew{\all n}}
M_{\Ew{\all n}}\herm V^{1/2}\bigr]\Bigr)}\norm{X(I-X)^{-1}}_2
+ \abs{\tr(\delta[A]_{(\Ew{\occ n, \occ n})})}\nonumber\\
&+ \frac12 \abs{\tr\Bigl(V^{1/2}M_{\Ew{\all n}}
M_{\Ew{\all n}}\herm V^{1/2}\Bigr)} \cdot 2\norm{\delta[X]}_2
+ o(\Fnorm{\delta[M]})\nonumber\\
\leq{}& \norm{\delta[X]}_2\tr\left({A_{(\all n, \all n)}}\right)
+ \abs{\tr(\delta[A]_{(\all n, \all n)})}\notag\\
&+ \frac{1}{2} \left\vert\tr[\delta[A]]\right\vert+o(\Fnorm{\delta[M]})\notag\\
={}& \norm{\delta[X]}_2\tr\left({A_{(\all n, \all n)}}\right)
+ \frac{3}{2} \abs{\tr(\delta[A_{(\all n, \all n)}])}
+ o(\Fnorm{\delta[M]}).\label{eq:rel_err_bound}
\end{align*}

\newpage

\footnotesize{
\bibliography{ref}
}